\documentclass[11pt]{article}
\usepackage{natbib}
\usepackage{geometry}                
\usepackage{graphicx}
\usepackage{amssymb}
\usepackage{epstopdf}
\title{Hierarchical Economic Agents and their Interactions}
\author{Ted Theodosopoulos\footnote{Department of Mathematics, Saint Ann's school, Brooklyn, NY, USA}}
\date{}                                           

\begin{document}
\maketitle

\begin{abstract}
We present a new type of spin market model, populated by hierarchical agents, represented as configurations of sites and arcs in an evolving network.  We describe two analytic techniques for investigating the asymptotic behavior of this model: one based on the spectral theory of Markov chains and another exploiting contingent submartingales to construct a deterministic cellular automaton that approximates the stochastic dynamics.  Our study of this system documents a phase transition between a sub-critical and a super-critical regime based on the values of a coupling constant that modulates the tradeoff between local majority and global minority forces.  In conclusion, we offer a speculative socioeconomic interpretation of the resulting distributional properties of the system.
\end{abstract}

\section{Introduction}
The primary goal of this paper is to describe the potential of a new economic modelling environment, populated by multi-layered agents, hierarchical objects that probe the boundary between individual and group, institution and society.  Bypassing questions of aggregation, the proposed paradigm seeks coordination through hierarchical, heterogeneous agents, influencing one-another through their opinions and actions.  Importantly, the agents' limited rationality permits pockets of inconsistent allegiances to percolate through their interaction network.

The proposed modeling environment extends work over the past decade on agent-based models of the economy [21,20,8,17].  Progressively, such models have shown how heterogeneities in the agents' endowments, preferences and interactions can persist and lead to observable deviations from the {\it efficient market hypothesis}, a collection of so-called {\it stylized facts} [16,13,18,7,26,2].  The extension proposed here invites us to broaden our notion of heterogeneity to encompass attributes that aren't reducible to individuals, but instead arise at different levels of aggregation.  But instead of seeking to extract them from properties of the individuals, we posit them as part of the evolving state of `meta-individualist', multi-layer agents that populate the economy [19,31].

To illustrate this broader modelling paradigm, we proceed to extend a specific agent-based model of the economy, first proposed by Bornholdt and subsequently studied both numerically and analytically by different authors [22,32].  This model is based on an interaction potential that trades off two components, the desire to belong to a local majority and simultaneously to the global minority.  These two terms are balanced by a {\it coupling constant} $\alpha$.  The study of the statistical mechanics of this model led to the identification of an explicit phase transition [32], controlled by $\alpha$, whereby sufficiently strong coupling leads to non-self-averaging behavior [3,4] and persistent opinion mixing.  Furthermore, this framework allowed us to identify a fundamental, irreducible limit to the observability of various measures of excess demand [32].

In the original Bornholdt model, states of the economy were represented by spin configurations of a fixed lattice, or network more generally.  Configurations of this kind can be denoted by vectors of $-1$s and $+1$s, $\eta \in \{-1,+1\}^N$, for some $N$.  These vectors are then propagated following a Markov process, driven by an interaction potential, $h(\eta)$.  The resulting stochastic evolution is seeking minimum energy states, which represent equilibria.  This minimization is controlled by a `temperature' parameter, in an analogy to the simulated annealing process of non-equilibrium statistical mechanics.  At high temperatures, state transitions are largely random.  As the temperature is lowered, transitions that locally reduce the interaction potential are progressively favored.  In the {\it frozen} phase, the system picks out some equilibrium state, and is subsequently trapped there.  

This framework is often interpreted as describing the evolution of individual agents, represented by the different sites on the lattice or network, with the spins at each site denoting the evolving opinions or actions (buys vs. sells) of the agent on that site.  However, such an interpretation, which attempts to reduce the resulting market dynamics to the interactions of individual agents, has come repeatedly under fire, from different perspectives [23,19].  Most problematic, from the point of view of the current paper, is the inability of this framework to reproduce any of the myriad intermediate structures, from coalitions to firms, that populate the real economic landscape.  More precisely, present instantiations of this paradigm reserve heterogeneity for individual agents, relegating any higher structures to the realm of transient epiphenomena.

Here, we propose to extend the standard spin market framework in an explicit effort to bring out the irreducible relevance of structures in the economy.  We choose to see these intermediate structures as endowed properties of the economic state, largely indecomposable to their constituents, albeit spontaneously evolving, in interaction among themselves and their constituents.  To help visualize our proposed scheme, we propose the following abstraction: agents are analogous to simplicial complexes in topology, consisting of locally matching components of different dimensionality and degree of complexity [27].  Such an object is generally indecomposable to a listing of its constituents.  Instead, it depends crucially on details of the `gluing map' that put it together.  Extending the analogy further, we posit a generalized interaction potential, which allows such hierarchical objects to `act' on one-another, without this action being describable as the interaction between individual components, e.g. an edge interacting with an edge or a tetrahedron interacting with a triangle.

The goal of this work is to introduce this new modeling paradigm, in which the agents and the network on which their opinions evolve are indissolubly coupled.  Unlike earlier agent-based studies of economic interactions, we don't attempt to generate a price process that can be calibrated against empirical statistics.  The translation of agents' opinions to observable aggregates depends sensitively on market microstructure, from details of the double-auction to explicit market making [9,31,1,6].  Instead, the current work focuses on the rich array of stochastic convergence effects that arise in models of heterogeneous economic agent, particularly when we endogenies the evolution of the interaction network [15,11].  Specifically, while much emphasis has been placed on conditions for guaranteeing ergodicity [5,28], the focus here is on the effective lack of ergodicity under certain parametric regimes for our model, and its economic consequences, both theoretical and empirical.  Along the way, we introduce techniques from symbolic dynamics and the spectral analysis of Markov chains to enrich the economic toolkit.  

We begin with an introduction to our hierarchical agent model, including a description of the interaction potential that couples the spin configurations on the nodes and arcs of our evolving network.  This interaction potential drives a Markov process whose hypergeometric state transitions are described in Section 3, along with some sample paths that hint at the non-ergodic behavior we are after.  The following section introduces the contingent submartingale representation, a technique that allows us to extract a deterministic skeleton underlying our stochastic dynamics.  We then proceed to investigate the invariant measure of the Markov process and its sensitivity to the model's parameters.  The discrepancies between the limiting distribution and the deterministic attractors we identified earlier leads us to pursue a spectral analysis of the underlying Markov chain, which uncovers and quantifies a source of persistent path dependence.  The paper concludes with a set of phenomenological conjectures that govern the paths of our hierarchical agent model, as well as a putative socioeconomic interpretation [10] for the three distinct dynamic regimes that the model exhibits.  All along, we relegate the more technically demanding details of our exposition to an Appendix.

\section{Model Description}
\begin{figure}[htbp] 
   \centering
   \includegraphics[width=4in]{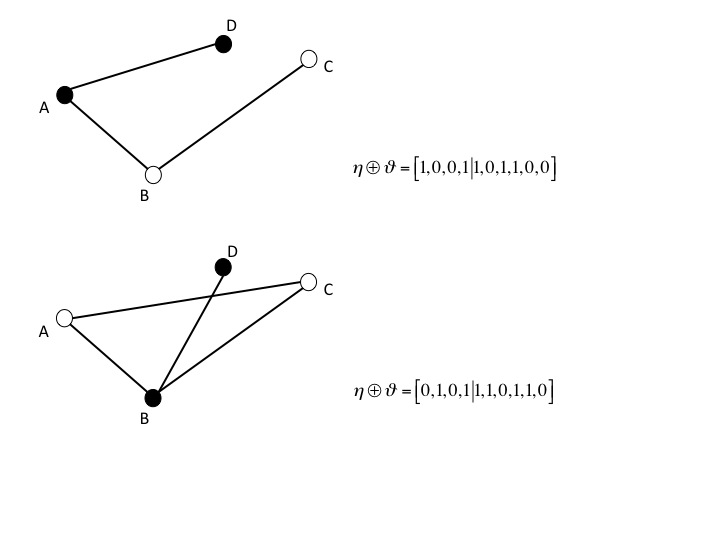} 
   \caption{Two different configurations of sites and arcs with $N = 4$.}
   \label{fig:examples}
\end{figure}

More concretely, we proceed to describe in detail an extension of the earlier Bornholdt spin market model, where the states of the economy are represented by an object with two components: binary configurations on sites and arcs.  To begin with, there is a set of sites ${\mathcal S} = \{1,2,\ldots,N\}$ for some $N$ and a related set of arcs ${\mathcal A} = \left\{ \left. (i,j) \in {\mathcal S}^2 \right| i<j \right\}$.  The first component of the state of the economy are spin configurations of sites, i.e. vectors of $-1$s and $+1$s, $\eta : {\mathcal S} \times (0,\infty) \rightarrow \{-1,+1\}$, while the second component are spin configurations over the set of arcs, $\vartheta: {\mathcal A} \times (0,\infty) \rightarrow \{-1,+1\}$.  Thus, the state can be described as 
$${\bf a} = \eta \oplus \vartheta = \left[\eta(1), \eta(2), \ldots, \eta(N) \left| \vartheta(1,2), \vartheta(1,3), \ldots, \vartheta(N-1,N) \right. \right],$$ 
as shown in the examples illustrated in Figure \ref{fig:examples}.  Let $S^+ (t) = \left| \left\{x \in {\mathcal S} \left| \eta(x,t) = +1 \right. \right\} \right|$ and $A^+ (t) = \left| \left\{(y,z) \in {\mathcal A} \left| \vartheta(y,z,t) = +1 \right. \right\} \right|$ denote the number of positive sites and arcs respectively.  For notational convenience, we extend $\vartheta$ to the whole ${\mathcal S}^2$, by imposing symmetry, i.e. $\vartheta (x,y,t) = \vartheta (y,x,t)$ for all $t$.

We construct a continuous time Markov process with transitions occurring at exponentially distributed epochs, $T_n$, with rate $1$ [14].  At time $T_n$ (i.e. the $n^{\rm th}$ epoch) a random member of ${\mathcal S} \times {\mathcal A}$ is chosen uniformly and its spin is changed to $+1$ or $-1$ depending on interactions between the two components of the current market configuration.  These interactions between these two components of the objects that populate our model rely on a tradeoff between the desire to align with the majority within a local neighborhood and a need to react to the opportunities created by global imbalance.  The neighborhood structure of sites is based on the current configuration of the arcs, and vice versa, exploiting the duality between sites and arcs.  In fact, as we will discuss in more detail later, this is one of the special features of our two-tiered agent that we explicitly exploit, and which is substantially more convoluted in higher-order extensions of our framework.  

More specifically, let ${\mathcal N}_\vartheta : {\mathcal S} \times (0,\infty) \rightarrow {\mathcal P}\left({\mathcal S}\right)$ denote a mapping that assigns to every site $y \in {\mathcal S}$ and every epoch $T_n$ a subset ${\mathcal N}_\vartheta (y,T_n)$ of ${\mathcal S}$ given by 
$$\left\{ \left. x \in {\mathcal S} \setminus \{y \} \right| \vartheta \left(x,y,T_n \right) = +1 \right\}.$$  
Similarly, let ${\mathcal N}_\eta : {\mathcal A} \times (0,\infty) \rightarrow {\mathcal P}({\mathcal A})$ denote a mapping that assigns to every arc $(x,y) \in {\mathcal A}$ and every epoch $T_n$ a subset ${\mathcal N}_\eta \left(x,y,T_n \right)$ of ${\mathcal A}$ given by ${_xA^y} \cup {^yA_x} \cup {_yA^x} \cup {^xA_y}$, where
\begin{eqnarray*}
_iA^j & = & \left\{ \left. (i,k) \in {\mathcal A} \right| k \ne j \mbox{ and } \eta (i) = +1 \right\} \\
^iA_j & = & \left\{ \left. (k,j) \in {\mathcal A} \right| k \ne i \mbox{ and } \eta (j) = +1 \right\}.
\end{eqnarray*}

The treatment of the model presented in this paper is based on applying a version of `rapid stirring' [12] by randomizing the neighborhood structure generated by ${\mathcal N}_\vartheta$ and ${\mathcal N}_\eta$.  Appendix 1 provides more details about how these random neighborhoods are drawn.  Once the neighborhoods have been assigned for each site and arc at a point in time, each member of each neighborhood is mapped to $+1$ or $-1$ using further hypergeometric random variables, independent of the earlier ones and of each other, with as many draws as there are members of the neighborhood under consideration.  These draws are without replacement, out of a population which depends on whether the chosen element is a site or an arc, and number of successes depending on the sign of the base site or arc.  In particular, the interaction potential for site $x \in {\mathcal S}$ and arc $(x,y) \in {\mathcal A}$ is given by
\begin{eqnarray}
h_{\eta \oplus \vartheta} (x,T_n) & = & \sum_{y \in {\cal N}_\vartheta(x,T_n)} \eta(y,T_n) - \alpha \eta(x,T_n) G(T_n) \nonumber \\
g_{\eta \oplus \vartheta} (x,y,T_n) & = & \sum_{(u,v) \in {\cal N}_\eta \left(x,y,T_n \right)} \vartheta(u,v,T_n) - \alpha \vartheta(x,y,T_n) G(T_n), \label{eq:potentials}
\end{eqnarray}
where
$$G(T_n) = {\frac {1}{2}} \left| {\frac {4 \left(S^+(T_n) + A^+(T_n) \right)}{N(N+1)}} -1 \right|$$
measures the global imbalance, in sites and arcs.

The dynamics of the state proceed as follows.  We set a temperature parameter, which controls the amount of randomization that interferes with the minimization of the interaction potentials described in (\ref{eq:potentials}) above.  As usual, we denote by $\beta$ the inverse temperature, and eventually we let it increase towards infinity.  At every point in time $n$, a site or arc is chosen uniformly at random and a coin is flipped.  The chosen site or arc is assigned a $+1$ if the coin comes up HEADS and $-1$ otherwise.  If site $x \in {\mathcal S}$ is chosen, then the probability that the coin comes up HEADS is equal to
$$\left(1+ \exp\left\{ -2\beta h_{\eta \oplus \vartheta} (x,T_n) \right\} \right)^{-1}.$$
On the other hand, if arc $(x,y)$ is chosen, then the probability that the coin comes up HEADS is equal to
$$\left(1+ \exp\left\{ -2\beta g_{\eta \oplus \vartheta} (x,y,T_n) \right\} \right)^{-1}.$$
Note that, naturally, arcs are chosen more often than sites, because there are quadratically more arcs than sites, but in the long run, this imbalance in refresh rates guarantees that the more arcs will have had an equal opportunity of settling down to their invariant marginal distribution as the significantly fewer sites.

\section{Transition Probabilities}
In order to proceed with our analysis, we will compute the probabilities that at any point in time, the site or arc that is chosen will not change its sign.  We will restrict our attention to the `frozen phase', i.e. we will consider the limit $\beta \rightarrow \infty$.  This choice simplifies our analysis because, in this limit, all `thermal' randomness, which would oppose the minimization of the interaction potential, disappears, and the only randomness that remains stems from the sampling of random neighborhoods.  This persistent randomization gives rise to hypergeometric random variables. 

More precisely, consider $P_{++}(i,j)$, the probability that, having chosen a $+1$ site, it remains $+1$ after the update, assuming the system is in a state with $S^+ = i$ and $A^+ = j$.  Let 
$$L(i,j,\ell) = 1+ \left \lfloor {\frac {1}{2}} \left( \ell+ {\frac {\alpha}{2}} \left|{\frac {4(i+j)}{N(N+1)}} -1 \right| \right) \right \rfloor,$$ 
where $\lfloor x \rfloor$ is the {\it floor} of $x$, i.e. the largest integer no greater than $x$.  Then, the site transition probability $P_{++}(i,j)$ is given by the following partial sum of conditionally hypergeometric random variables:
\begin{equation}
P_{++} (i,j) = \left(\begin{array}{c} {\rm C}^N_2 \\ N-1 \end{array} \right)^{-1} \sum_{\ell = 0 \vee (j- {\rm C}_2^{N-1})}^{j \wedge (N-1)} a_\ell \sum_{k = L(i,j,\ell)}^{\ell \wedge (i-1)} b_k, \label{eq:ppp}
\end{equation}
where 
$$a_\ell = \left(\begin{array}{c} j \\ \ell \end{array} \right) \left(\begin{array}{c} {\rm C}^N_2 -j \\ N-\ell-1 \end{array} \right) \left(\begin{array}{c} N-1 \\ \ell \end{array} \right)^{-1},$$ 
$$b_k = \left(\begin{array}{c} i-1 \\ k \end{array} \right) \left(\begin{array}{c} N-i \\ \ell-k \end{array} \right)$$
and $a \vee b = \max\{a,b\}$ and $a \wedge b = \min\{a,b\}$.  More details about the derivation of these transition probabilities are given in Appendix 2. 

Similarly the probability $P_{--} (i,j)$ that, having chosen a $-1$ site, it remains $-1$ after the update, assuming the system is in a state with $S^+ = i$ and $A^+ = j$, is based on guaranteeing that $h_{\eta \oplus \vartheta} (x)<0$ and is given by
\begin{equation}
P_{--} (i,j) = \left(\begin{array}{c} {\rm C}^N_2 \\ N-1 \end{array} \right)^{-1} \sum_{\ell = 0 \vee (j- {\rm C}_2^{N-1})}^{j \wedge (N-1)} c_\ell \sum_{k = 0 \vee (i+\ell+1-N)}^{U(i,j,\ell)} d_k, \label{eq:pmm}
\end{equation}
where $U(i,j,\ell) = \left \lceil {\frac {1}{2}} \left( \ell- {\frac {\alpha}{2}} \left|{\frac {4(i+j)}{N(N+1)}} -1 \right| \right) \right \rceil-1$, 
$$c_\ell = \left(\begin{array}{c} j \\ \ell \end{array} \right) \left(\begin{array}{c} {\rm C}^N_2 -j \\ N-\ell-1 \end{array} \right) \left(\begin{array}{c} N-1 \\ \ell \end{array} \right)^{-1},$$
$$d_k = \left(\begin{array}{c} i \\ k \end{array} \right) \left(\begin{array}{c} N-i-1 \\ \ell-k \end{array} \right)$$
and $\lceil x \rceil$ is the {\it ceiling} of $x$, i.e. the smallest integer no less than $x$.

On the other hand, if we choose an arc instead of a site, the computation is somewhat different.  As before, we are looking to compute the probability, $Q_{++} (i,j)$ that, having chosen a $+1$ arc, it remains $+1$ after the update, assuming the system is in a state with $S^+ = i$ and $A^+ = j$.  Let
$$R(i,j) = 1+\left \lfloor {\frac {1}{2}} \left(N-2+ {\frac {\alpha}{2}} \left|{\frac {4(i+j)}{N(N+1)}} -1 \right| \right) \right \rfloor$$
and
$$S(i,j) = 1+\left \lfloor {\frac {1}{2}} \left(2N-4+ {\frac {\alpha}{2}} \left|{\frac {4(i+j)}{N(N+1)}} -1 \right| \right) \right \rfloor.$$
Then, the arc transition probability $Q_{++}(i,j)$ is given by the following partial sum of conditionally hypergeometric random variables:
\begin{eqnarray}
Q_{++} (i,j) & = & {\frac {2i(N-i)}{N(N-1)}} \left(\begin{array}{c} {\rm C}^N_2 -1 \\ N-2 \end{array} \right)^{-1} \sum_{k=R(i,j) \vee \left(j - {\frac {N^2 - 3N + 4}{2}} \right)}^{(j-1) \wedge (N-2)} \left(\begin{array}{c} j-1 \\ k \end{array} \right) \left(\begin{array}{c} {\rm C}^N_2 -j \\ N-k-2 \end{array} \right) + \nonumber \\
& & + {\frac {i(i-1)}{N(N-1)}} \left(\begin{array}{c} {\rm C}^N_2 -1 \\ 2N-4 \end{array} \right)^{-1} \sum_{k=S(i,j) \vee \left(j - {\frac {N^2 - 5N + 8}{2}} \right)}^{(j-1) \wedge (2N-4)} \left(\begin{array}{c} j-1 \\ k \end{array} \right) \left(\begin{array}{c} {\rm C}^N_2 -j \\ 2N-k-4 \end{array} \right). \nonumber \\
& & \label{eq:qpp}
\end{eqnarray}

Similarly, if the chosen arc is $-1$, the probability that it remains negative after the update is given by 
\begin{eqnarray}
Q_{--} (i,j) & = & {\frac {2i(N-i)}{N(N-1)}} \left(\begin{array}{c} {\rm C}^N_2 -1 \\ N-2 \end{array} \right)^{-1} \sum_{k=0 \vee \left(j- {\frac {N^2-3N+2}{2}} \right)}^{T(i,j)} \left(\begin{array}{c} j \\ k \end{array} \right) \left(\begin{array}{c} {\rm C}^N_2 -j-1 \\ N-k-2 \end{array} \right) + \nonumber \\
& & + {\frac {i(i-1)}{N(N-1)}} \left(\begin{array}{c} {\rm C}^N_2 -1 \\ 2N-4 \end{array} \right)^{-1} \sum_{k=0 \vee \left(j- {\frac {N^2-5N+6}{2}} \right)}^{V(i,j)} \left(\begin{array}{c} j \\ k \end{array} \right) \left(\begin{array}{c} {\rm C}^N_2 -j-1 \\ 2N-k-4 \end{array} \right), \nonumber \\
& & \label{eq:qmm}
\end{eqnarray}
where
$$T(i,j) = \left \lceil {\frac {1}{2}} \left(N-2- {\frac {\alpha}{2}} \left|{\frac {4(i+j)}{N(N+1)}} -1 \right| \right) \right \rceil -1$$
and
$$V(i,j) = \left \lceil {\frac {1}{2}} \left(2N-4- {\frac {\alpha}{2}} \left|{\frac {4(i+j)}{N(N+1)}} -1 \right| \right) \right \rceil -1.$$
The technical details of these derivations are also described in Appendix 2. 

\begin{figure}[htbp] 
   \centering
   \includegraphics[width=5in]{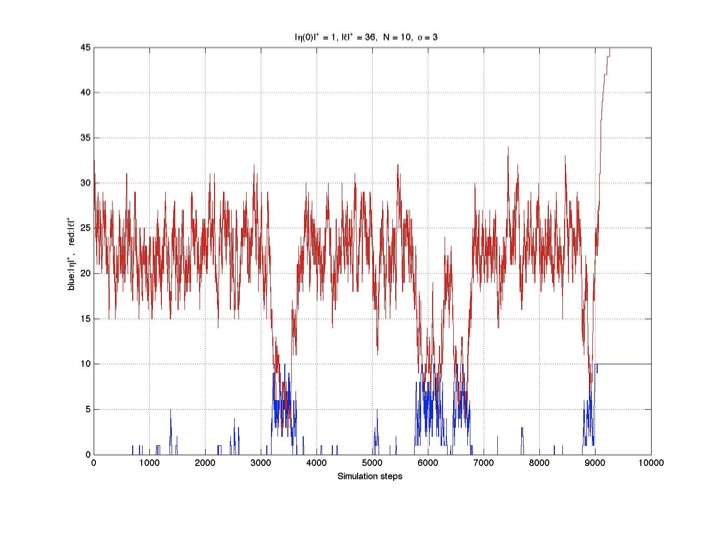} 
   \caption{Monte Carlo simulation of 2D Markov Chain with $N=10$ and $\alpha = 3$.}
   \label{fig:MC1}
\end{figure}

We can perform a Monte Carlo simulation based on the expressions (\ref{eq:ppp}), (\ref{eq:pmm}), (\ref{eq:qpp}) and (\ref{eq:qmm}).  In particular, we can consider a two-dimensional Markov Chain ${\bf X}_n = \left(S^+(T_n), A^+(T_n) \right)$ on $\{1,2,\ldots,N\} \times \left\{1,2,\ldots,{\rm C}^N_2 \right\}$, with transition probabilities given by
\begin{eqnarray}
& & {\rm Pr} \left({\bf X}_{n+1} = (k,\ell) \left| {\bf X}_n = (i,j) \right. \right) = \left\{ \begin{array}{cc} 
{\frac {2(N-i)}{N(N+1)}}\left[1-P_{--} (i,j) \right] & \mbox{if $k = i+1$ and $\ell = j$} \\
{\frac {2i}{N(N+1)}}\left[1-P_{++} (i,j) \right] & \mbox{if $k = i-1$ and $\ell = j$} \\
{\frac {N^2-N-2j}{N(N+1)}}\left[1-Q_{--} (i,j) \right] & \mbox{if $k = i$ and $\ell = j+1$} \\
{\frac {2j}{N(N+1)}}\left[1-Q_{++} (i,j) \right] & \mbox{if $k = i$ and $\ell = j-1$} \\
{\frac {2L(i,j)}{N(N+1)}} & \mbox{if $k = i$ and $\ell = j$}
\end{array} \right. \nonumber \\
& & \label{eq:MCtransitions}
\end{eqnarray}
where 
$$L(i,j) = NP_{--}(i,j) + C_2^NQ_{--}(i,j) + i\left[ P_{++}(i,j) - P_{--}(i,j) \right] + j\left[ Q_{++}(i,j) - Q_{--}(i,j) \right].$$
Figure \ref{fig:MC1} shows a simulation of this Markov Chain for $10,000$ steps, with $N=10$ and $\alpha = 3$, starting at $S^+(0) = 1$ and $A^+(0) = 36$.  After about $9,000$ steps, this simulation was trapped in the attractor ${\bf X} = (10,45)$.

\begin{figure}[htbp] 
   \centering
   \includegraphics[width=5in]{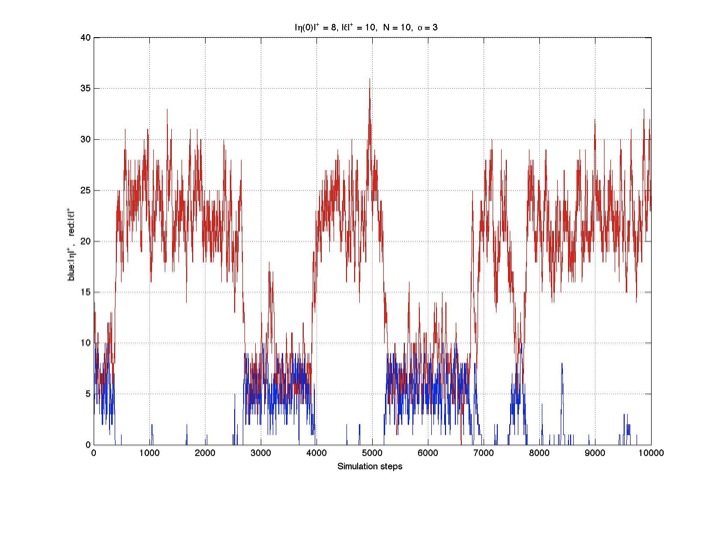} 
   \caption{Another Monte Carlo simulation of 2D Markov Chain with $N=10$ and $\alpha = 3$, started at a different point.}
   \label{fig:MC2}
\end{figure}

This Markov Chain is a form of two-dimensional random walk in a random environment (RWRE), i.e. a random walk in which the probability of going UP/DOWN and LEFT/RIGHT depends on your current location [12].  Such stochastic processes may exhibit path dependence and lack of ergodicity, as supported by the simulation of the same Markov Chain shown in figure \ref{fig:MC2}.

This simulation also lasts for $10,000$ steps, and the parameters are identical with the earlier simulation.  The only difference is that this time the simulation was started at $S^+(0) = 8$ and $A^+(0) = 10$.  Instead of becoming trapped in the attractor ${\bf X} = (10,45)$, this time the Markov chain appears to be stochastically switching between two states, one with $(S^+,A^+) \approx (5,8)$ and another with $(S^+,A^+) \approx (0,22)$.  This phenomenology is consistent with a non-ergodic process, for which the dependence on initial conditions doesn't disappear in the $n \rightarrow \infty$ limit.

One of the most prominent characteristics of empirical price series is their long range memory [9,1].  This can indicate path-dependence, a tell-tale sign of non-ergodicity in the underlying process.  In lieu of generating price dynamics, which could be compared with empirically determined statistics, we explore the serial correlation of the $S^+$ and $A^+$ paths generated by our process for signs of similar qualitative behavior.  We quantify the serial correlation using the modified R/S index [25,24].  This index is a function of a delay window, $\tau$, and is normalized so that when its value is equal to $0.5$ for a particular delay, the process is memoryless at that horizon.  On the other hand, when the R/S index is more than $0.5$, the process exhibits persistent behavior, while values of the index below $0.5$ indicate anti-persistent behavior, at the corresponding horizons.

\begin{figure}[htbp] 
   \centering
   \includegraphics[width=3.5in]{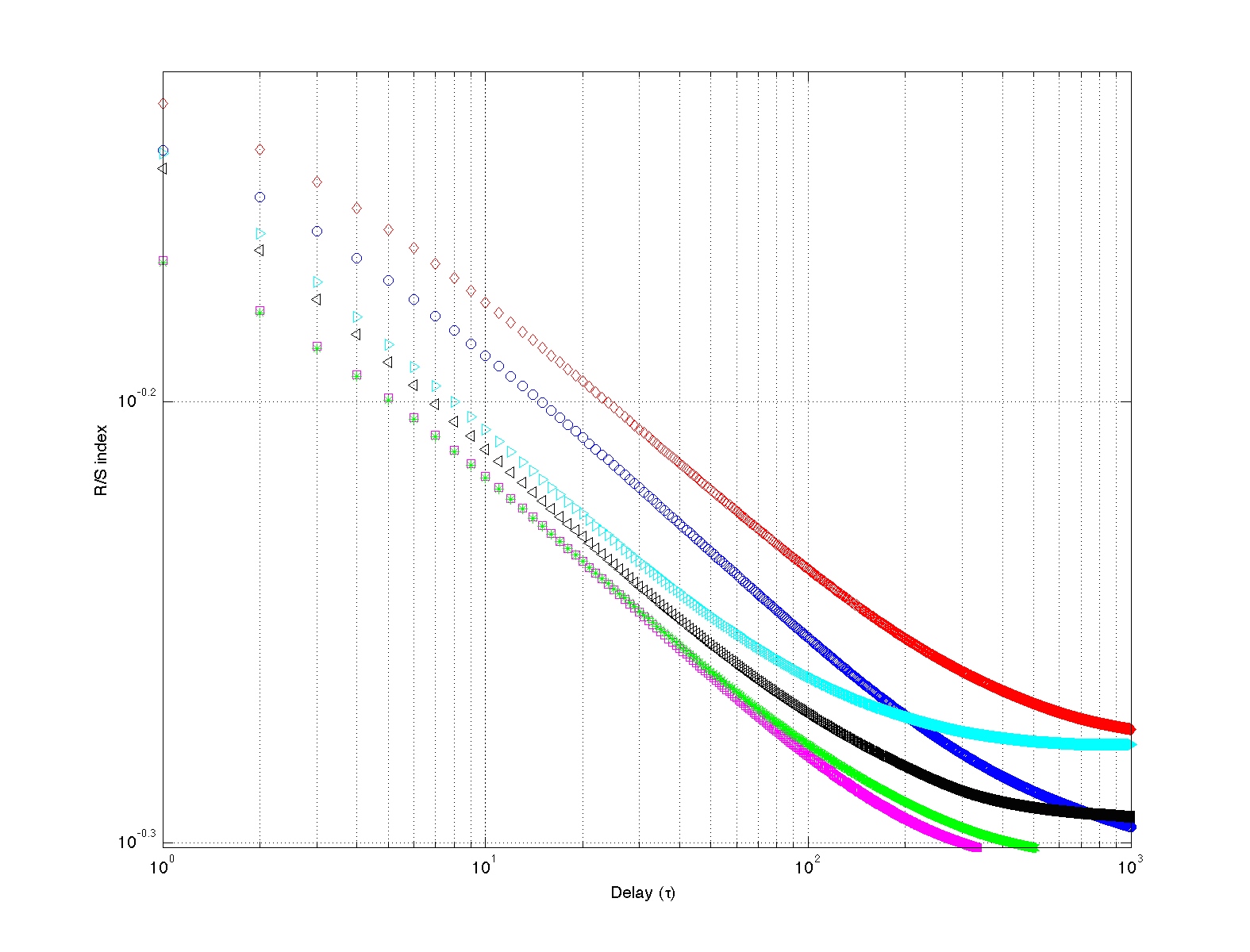} 
   \caption{The stochastic paths of $S^+$ and $A^+$ exhibit persistent serial correlation, as shown by the slowly decaying R/S index.  The lower edge corresponds to a value of R/S equal to $0.5$.  $S^+$: Blue Cirlces ($\alpha = 10$), Magenta Squares ($\alpha = 20$), Black Left Triangles ($\alpha = 40$; $A^+$: Red Diamonds ($\alpha = 10$), Green Stars ($\alpha = 20$), Cyan Right Triangles ($\alpha = 40$).}
   \label{fig:RS1}
\end{figure}

Figure \ref{fig:RS1} shows the R/S index for three Monte Carlo simulation runs, each with $100,000$ time steps.  Note that the figure is cropped so that the bottom edge corresponds to a value of $0.5$.  In all three cases $N = 10$ and the Markov Chain was started at $S^+(0) = 8$ and $A^+(0) = 10$.  Three different values of $\alpha$ were used: $10$, $20$ and $40$ respectively.  In the first and the last case, both $S^+$ and $A^+$ exhibit strong persistence even at $\tau>1,000$.  In the middle case, the memory disappears for both $S^+$ and $A^+$ after about $300$ and $500$ steps respectively.  It is worth noting that in all cases the arc process has longer memory than the site process.  

\begin{figure}[htbp] 
   \centering
   \includegraphics[width=5in]{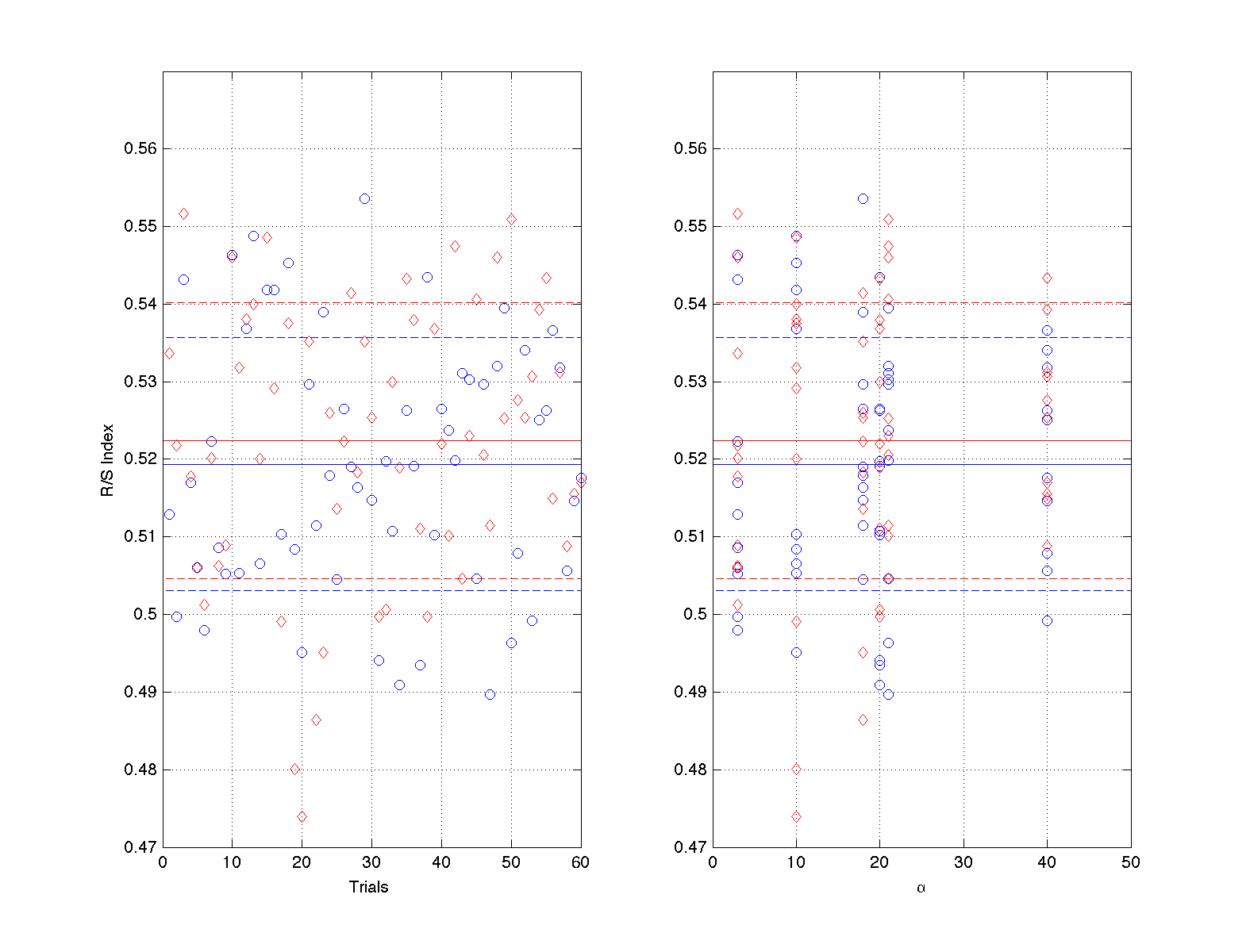} 
   \caption{Over 60 simulations with different values of the coupling constant $\alpha$, the R/S index remained above $0.5$ more than $88\%$ of the time for both $S^+$ (Blue Circles) and $A^+$ (Red Diamonds).  The left panel shows the resulting 120 values of the R/S index by trial, while the right panel shows the same data by value of the coupling constant $\alpha$.}
   \label{fig:RS2}
\end{figure}

To further explore the long range memory of our stochastic process, we ran $60$ Monte Carlo simulations with $10,000$ steps each, all started at $S^+(0) = 8$ and $A^+(0) = 10$, with six different values of $\alpha$.  Figure \ref{fig:RS2} shows the R/S index at $\tau = 1,000$.  In each case the solid lines represent the mean values of the index and the dashed lines correspond to one standard deviation above and below the mean (blue corresponds to $S^+$ and red to $A^+$).  More than $85\%$ of the data remain persistently serially correlated even beyond the $1,000$ step horizon, indicating a remarkably long memory, which may well contribute to the observed long memory of empirical economic price series.

In what follows, we will employ a novel methodology to probe the dynamical attractors of this stochastic process.  This technique, based on a {\it contingent submartingale representation} [33], will allow us to construct a {\it cellular automaton} [34] approximation of the full stochastic dynamics, identify all equilibria of this deterministic dynamical system and examine their stability properties.  For linear stochastic systems, the resulting deterministic paths act as a `skeleton' around which the stochastic state paths oscillate, ultimately converging to invariant measures supported in the neighborhoods of the deterministic attractors.  

On the other hand, stochastic systems with nonlinear interaction terms like those in our system often possess more complex, non-classical limiting behaviors in path space, e.g. involving path dependence and lack of ergodicity, that aren't reducible to limiting averages over longer times.  In the following sections of this paper we will illustrate qualitative deviations between the paths of states (primal objects), which exhibit a bewildering array of asymptotic behaviors, from limit point to intertwining periodic orbits and chaotic attractors, and the evolution of measures (dual objects), which unambiguously converge to well-described distributions, whose properties and convergence rates we are in a position to characterize.

\section{Contingent Submartingale Representation}
In order to analyze the long term behavior of this decidedly complicated, inhomogeneous Markov chain, we begin by characterizing the conditional expectation of the increments of this process in each dimension.  This analysis gives us the following deterministic nonlinear dynamics:
\begin{eqnarray}
f(i,j) & \stackrel{\Delta}{=} & \left({\frac {N+1}{2}}\right) {\bf E} \left[S^+_{n+1} - S^+_n \left| {\bf X}_n = (i,j) \right. \right] \nonumber \\
& = & \left(1 - {\frac {i}{N}} \right) \left( 1- P_{--}(i,j) \right) - {\frac {i}{N}} \left( 1- P_{++}(i,j) \right) \nonumber \\
\label{eq:etaincr} \\
g(i,j) & \stackrel{\Delta}{=} & \left({\frac {N+1}{N-1}}\right) {\bf E} \left[A^+_{n+1} - A^+_n \left| {\bf X}_n = (i,j) \right. \right] \nonumber \\
& = & \left(1 - {\frac {j}{{\rm C}^N_2}} \right) \left( 1- Q_{--}(i,j) \right) - {\frac {j}{{\rm C}^N_2}} \left( 1- Q_{++}(i,j) \right), \nonumber \\
\label{eq:thetaincr} 
\end{eqnarray}
where we've taken to using $S^+_{k} = S^+ \left(T_k \right)$ and $A^+_{k} = A^+ \left(T_k \right)$ for notational simplicity.  In general, we are interested in the sign of $f$ and $g$, because that determines whether the two components of ${\bf X}_n$ increase or decrease, on average.  

Incidentally, the fractional factors in front of the conditional expectations reflect the fact that the moves in ${\bf X}$ will be either in the first ($S^+$) or second ($A^+$) dimension, not both simultaneously, with proportions $N : {\rm C}^N_2$, as we mentioned already at the end of section 2.  In other words, the arc configurations are updated more frequently than the site configurations, and therefore the conditional expectations of the increments of $S^+$ will be much lower than the ones for $A^+$, and they both will be lower than they would have been if we allowed simultaneous moves in both directions.  It is these latter (simultaneous) expected increments that the functions $f$ and $g$ compute, so they need to be adjusted accordingly, to bring the weighted averages in line.  

The reason why we chose to define the functions $f$ and $g$ in this way, despite the superficial conflict with the definition of our 2D stochastic process, is because we intend them to banish all stochasticity and represent the deterministic kernel around which the stochastic process evolves, as discussed in the last paragraph of section 3.  Since both $f$ and $g$ depend on both components, they will generally point in directions that combine movement in both directions.  Had we insisted to choose one direction over the other at every step, we would have had to introduce some randomization scheme that chooses the directional signal from the pair of $f$ and $g$ to obey at every point in time.  The only way to extract a completely deterministic dynamic is to accept the possibility of simultaneous (though deterministic) moves in both directions.  

More specifically, it is situations like this that gave rise to the concept of a {\it contingent submartingale} [33].  Here we slightly generalize the definition given in [32] to accommodate our needs\footnote{Strictly speaking, there is no generalization involved.  The use of a 2D contingency criterion for determining the martingale properties of a 1D process was allowed even in the original scheme [32], even one involving this 1D process itself as one of the components of the 2D contingency criterion.  But we adjust our notation of the contingent submartingale here to make this possibility more explicit, because it is of direct relevance.}.  Specifically, let $X_n$ and $Y_n$ be two stochastic processes, and $A$ be a subset of the range of $\left(X_n,Y_n \right)$.  Let $R_0 = S_0 = 0$ and define for $k \geq 1$ the following two sequences of integer-valued stopping times:
\begin{eqnarray*}
R_k = \inf \{ n>S_{k-1} \left| \left(X_n, Y_n \right) \in A \right. \}, \\
S_k = \inf \{ n>R_k \left| \left(X_n, Y_n \right) \not \in A \right. \}.
\end{eqnarray*}
Finally, consider a two-dimensional process $W: {\mathbb Z}_+ \times {\mathbb Z}_+ \rightarrow {\mathbb R}$ defined by
\begin{eqnarray*}
W_{0,k} = Y_{R_k}, \\
W_{n,k} = Y_{(R_k+n) \wedge S_k}.
\end{eqnarray*}
We will say that $Y$ is a {\it contingent submartingale} with respect to $\left((X,Y),A \right)$ if $W_{\cdot,k}$ is a submartingale for each $k$ [14].

\begin{figure}[htbp] 
   \centering
   \includegraphics[width=6.5in]{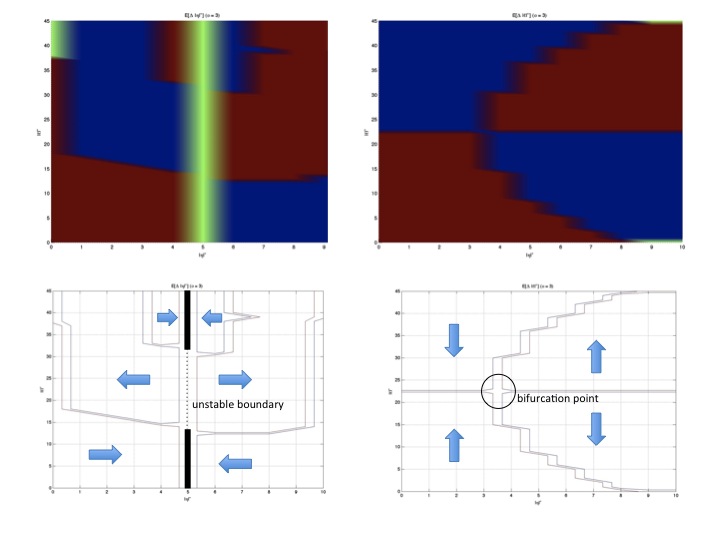} 
   \caption{Signs of the conditional first moments of the increments of ${\bf X}_n$ for $N = 10$ and $\alpha  = 3$.}
   \label{fig:CS1}
\end{figure}

Figure \ref{fig:CS1} shows the signs of the conditional expectations in (\ref{eq:etaincr}) and (\ref{eq:thetaincr}).  In particular, the top left panel shows the sign of the (conditionally expected) increments in $S^+$, with red, blue and yellow indicating positive, negative and zero respectively.  The top right panel similarly illustrates the sign of the (conditionally expected) increments in $A^+$.  In particular, if $U$ is the red region in the top left panel, then we can see that $S^+$ is a contingent submartingale with respect to $({\bf X}, U)$.  Similarly, if $V$ is the red region in the top right panel, then $A^+$ is a contingent submartingale with respect to $({\bf X}, V)$.

The bottom two panels indicate the boundaries between the regions of different signs.  Moreover, they include arrows to better illustrate the direction of the expected flows at each point in the state space.  Specifically, the boundaries of the different regions represent states where the process behaves locally like a martingale.  Intersections of the boundaries in the left and right bottom panels indicate states that are stationary, from an expectation point of view.  

\begin{figure}[htbp] 
   \centering
   \includegraphics[width=7in]{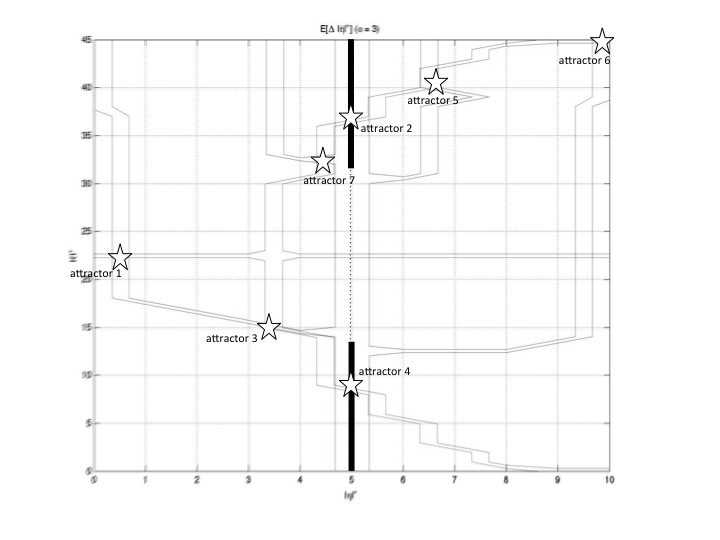} 
   \caption{Attractors of the cellular automaton abstraction.}
   \label{fig:CS2}
\end{figure}

At this point, it may be instructive to consider a thought experiment.  Imagine that we could do away with all residual randomness, and allow the system to follow purely the deterministic dynamics indicated by the vector fields in (\ref{eq:etaincr}) and (\ref{eq:thetaincr}) as shown in figure \ref{fig:CS1}:
\begin{eqnarray}
x_{n+1} = x_n + {\rm sgn} \left\{ f \left(x_n, y_n \right) \right\} \label{eq:xdeterm} \\
y_{n+1} = y_n + {\rm sgn} \left\{ g \left(x_n, y_n \right) \right\}, \label{eq:ydeterm} 
\end{eqnarray}
where ${\rm sgn}(a) = \left\{ \begin{array}{ll} a/|a| & \mbox{if $a \ne 0$} \\
0 & \mbox{otherwise} \end{array} \right.$ is the {\it signum} function.  This discrete dynamical system is what mathematicians call a {\it cellular automaton} [34].  This cellular automaton possesses several fixed point equilibria, located at the intersections of the curves $f = 0$ and $g = 0$.  Computationally, we can specify the approximate grid locations of the seven such equilibria when $\alpha = 3$
\begin{eqnarray*}
\left(x^\ast_1, y^\ast_1 \right) & = & (0, 23) \\
\left(x^\ast_2, y^\ast_2 \right) & = & (5, 36) \\
\left(x^\ast_3, y^\ast_3 \right) & = & (4, 15) \\
\left(x^\ast_4, y^\ast_4 \right) & = & (5, 8) \\
\left(x^\ast_5, y^\ast_5 \right) & = & (6, 41) \\
\left(x^\ast_6, y^\ast_6 \right) & = & (10, 45) \\
\left(x^\ast_7, y^\ast_7 \right) & = & (4, 33), 
\end{eqnarray*}
as can be seen in figure \ref{fig:CS2}, which superimposes the two lower panels in figure \ref{fig:CS1} and indicates their seven intersection points.

The question we now have to ask ourselves is whether these attractors are really there (because of the discretization effects which may diffuse the resulting local minima) and what their stability properties are.  More specifically, there are two separate effects to consider.  The first is due to the discrete nature of the underlying state space, which we have disregarded in our contingent submartingale analysis.  It is conceivable that a non-lattice point equilibrium may be unstable in the context of the discrete dynamics because all neighboring lattice points are repelling.  The second effect we need to consider is whether small perturbations would irreversibly escape the neighborhood of the putative attractor. 

Table 1 shows the values of $f$ and $g$ at each of the putative attractors and at grid points in their neighborhood, and uses this information to assess their stability properties.  The table shows that there is only on stable point attractor, namely $(10,45)$, while there are also three separate stable period-2 attractors for the dynamics when $\alpha = 3$.

\begin{table}\label{table:table1}
\begin{tabular}{c|c|c} 
Attractor & Lattice neighbors & Stability assessment \\ \hline
$\begin{array}{c}
(0,23) \leftrightarrows (1,22) \\
(1,23) \leftrightarrows (0,22)
\end{array}$ &
$\begin{array}{c} 
(f(0,21), g(0,21)) = (0.001, 0.0667) \\
(f(0,22), g(0,22)) = (0.0009, 0.0222) \\
(f(0,23), g(0,23)) = (0.0005, -0.0222) \\
(f(0,24), g(0,24)) = (0.0003,- 0.0667) \\
(f(1,21), g(1,21)) = (-0.0805, 0.0468) \\
(f(1,22), g(1,22)) = (-0.0852, 0.0156) \\
(f(1,23), g(1,23)) = (-0.0889, -0.0156) \\
(f(1,24), g(1,24)) = (-0.0919, -0.0468) \\
(f(2,21), g(2,21)) = (-0.1163, 0.0281) \\
(f(2,22), g(2,22)) = (-0.1268, 0.0093) \\
(f(2,23), g(2,23)) = (-0.1366, -0.0093) \\
(f(2,24), g(2,24)) = (-0.1458, -0.0281) 
\end{array}$ & 
$\begin{array}{c}
\mbox{Two stable periodic attractors} \\
\mbox{(period 2)}
\end{array}$\\ \hline
(5,36) & 
$\begin{array}{c} 
(f(4,35), g(4,35)) = (-0.1493, -0.0644) \\
(f(4,36), g(4,36)) = (-0.1574, -0.085) \\
(f(4,37), g(4,37)) = (-0.1663, -0.1076) \\
(f(5,35), g(5,35)) = (0, 0.0207) \\
(f(5,37), g(5,37)) = (0, -0.0167) \\
(f(6,35), g(6,35)) = (0.1493, 0.0908) \\
(f(6,36), g(6,36)) = (0.1574, 0.075) \\
(f(6,37), g(6,37)) = (0.1663, 0.0568) \end{array}$ & 
$\begin{array}{c}
\mbox{Saddle point } \\
\mbox{(unstable along horizontal axis)}
\end{array}$\\ \hline
(4,15) & 
$\begin{array}{c} 
(f(3,14), g(3,14)) = (-0.0159, 0.0854) \\
(f(3,15), g(3,15)) = (-0.0288, 0.0704) \\
(f(3,16), g(3,16)) = (-0.0414, 0.0572) \\
(f(4,14), g(4,14)) = (-0.0076, 0.0065) \\
(f(4,15), g(4,15)) = (-0.0145, -0.0016) \\
(f(4,16), g(4,16)) = (-0.0213, -0.0073) \\
(f(5,14), g(5,14)) = (0, -0.0632) \\
(f(5,15), g(5,15)) = (0, -0.0659) \\
(f(5,16), g(5,16)) = (0, -0.0653) 
\end{array}$ & 
$\begin{array}{c}
\mbox{Unstable} \\
\mbox{(escape through $(3,14)$)}
\end{array}$\\ \hline
$(5,8) \leftrightarrows (5,9)$ &
$\begin{array}{c} 
(f(4,7), g(4,7)) = (0.0493, 0.1323) \\
(f(4,8), g(4,8)) = (0.0395, 0.1076) \\
(f(4,9), g(4,9)) = (0.0306, 0.085) \\
(f(5,7), g(5,7)) = (0, 0.0387) \\
(f(5,8), g(5,8)) = (0, 0.0167) \\
(f(5,9), g(5,9)) = (0, -0.0031) \\
(f(5,10), g(5,10)) = (0, -0.0207) \\
(f(6,8), g(6,8)) = (-0.0395, -0.0568) \\
(f(6,9), g(6,9)) = (-0.0306, -0.075) \\
(f(6,10), g(6,10)) = (-0.0223, -0.0908)
\end{array}$ & 
Stable periodic attractor (period 2)
\\ \hline
\end{tabular}
\end{table}

\begin{table}\label{table:table2}
\begin{tabular}{c|c|c} 
Attractor & Lattice neighbors & Stability assessment \\ \hline
(6,41) & 
$\begin{array}{c} 
(f(5,40), g(5,40)) = (0, -0.0876) \\
(f(5,41), g(5,41)) = (0, -0.1138) \\
(f(5,42), g(5,42)) = (0, -0.1407) \\
(f(6,40), g(6,40)) = (-0.0523, -0.0085) \\
(f(6,41), g(6,41)) = (-0.0648, -0.0328) \\
(f(6,42), g(6,42)) = (-0.0829, -0.0578) \\
(f(7,40), g(7,40)) = (0.2298, 0.0509) \\
(f(7,41), g(7,41)) = (0.2519, 0.0279) \\
(f(7,42), g(7,42)) = (0.2724, 0.0044) \end{array}$ & 
$\begin{array}{c}
\mbox{Unstable} \\
\mbox{(escape from everywhere)}
\end{array}$\\ \hline
(10,45) & 
$\begin{array}{c} 
(f(9,44), g(9,44)) = (0.1, 0.0222) \\
(f(9,45), g(9,45)) = (0.1, 0) \\
(f(10,44), g(10,44)) = (0, 0.0222) \end{array}$ & 
Stable\\ \hline
(4,33) & 
$\begin{array}{c} 
(f(3,32), g(3,32)) = (-0.2207, -0.1023) \\
(f(3,33), g(3,33)) = (-0.2315, -0.1211) \\
(f(3,34), g(3,34)) = (-0.2421, -0.1419) \\
(f(4,32), g(4,32)) = (-0.1274, -0.0172) \\
(f(4,33), g(4,33)) = (-0.1344, -0.0305) \\
(f(4,34), g(4,34)) = (-0.1417, -0.0462) \\
(f(5,32), g(5,32)) = (0, 0.0572) \\
(f(5,33), g(5,33)) = (0, 0.0479) \\
(f(5,34), g(5,34)) = (0, 0.0357) 
\end{array}$ & 
$\begin{array}{c}
\mbox{Unstable} \\
\mbox{(escape from everywhere)}
\end{array}$\\ \hline
\end{tabular}
\caption{Stability analysis of putative attractors of the cellular automaton model (\ref{eq:xdeterm}) and (\ref{eq:ydeterm})}
\end{table}

So, the deterministic dynamics of the cellular automaton approximation are non-ergodic, because they lack uniqueness of asymptotic behavior.  The effect of initial conditions never disappears.  Next we proceed to shift our perspective from the time evolution of the individual paths the system may follow to the evolution of measures on the space ${\mathcal S} \times {\mathcal A}$.  We investigate the fully stochastic system by mapping it into a Markov chain and using spectral methods to study its convergence properties.  Throughout the following discussion, it is instructive to keep in mind that the objects propagated by the Markov chain aren't individual states, but measures [29].  

\section{Invariant Measure}
In order to take advantage of linear algebra, we will recast the 2D Markov chain in terms of propagating state vectors.  Specifically, let 
\begin{equation}
J(x) = \left \lfloor {\frac {x}{{\rm C}_2^N + 1}} \right \rfloor \label{eq:Jfunction}
\end{equation} 
and 
\begin{equation}
I(x) = x - J(x) \left({\rm C}^N_2 + 1 \right). \label{eq:Ifunction}
\end{equation}
Then, using (\ref{eq:MCtransitions}) we obtain
\begin{eqnarray}
& & M_{ab} = \left\{ \begin{array}{cc} 
{\frac {2(N-I(a))}{N(N+1)}}\left[1-P_{--} (I(a),J(a)) \right] & \mbox{if $I(b) = I(a) +1$ and $J(b) = J(a)$} \\
{\frac {2I(a)}{N(N+1)}}\left[1-P_{++} (I(a),J(a)) \right] & \mbox{if $I(b) = I(a)-1$ and $J(b) = J(a)$} \\
{\frac {N^2-N-2J(a)}{N(N+1)}}\left[1-Q_{--} (I(a),J(a)) \right] & \mbox{if $I(b) = I(a)$ and $J(b) = J(a)+1$} \\
{\frac {2J(a)}{N(N+1)}}\left[1-Q_{++} (I(a),J(a)) \right] & \mbox{if $I(b) = I(a)$ and $J(b) = J(a)-1$} \\
{\frac {2L(I(a),J(a))}{N(N+1)}} & \mbox{if $I(b) = I(a)$ and $J(b) = J(a)$}
\end{array} \right. \nonumber \\
& & \label{eq:Matrix}
\end{eqnarray}
where $L(i,j)$ is defined after (\ref{eq:MCtransitions}).  For example, consider the case $N = 3$ and $\alpha= 6$.  The state is described by elements of $\{0,1,2,3\} \times \{0,1,2,3\}$, a set with cardinality 16.  Thus, the resulting $16 \times 16$ state transition matrix $M$ will be given by:
\begin{equation}
M = \left( \begin{array}{cccccccccccccccc}
0 & {\frac {1}{2}} & 0 & 0 & {\frac {1}{2}} & 0 & 0 & 0 & 0 & 0 & 0 & 0 & 0 & 0 & 0 & 0 \\
{\frac {1}{6}} & 0 & {\frac {1}{3}} & 0 & 0 & {\frac {1}{2}} & 0 & 0 & 0 & 0 & 0 & 0 & 0 & 0 & 0 & 0 \\
0 & {\frac {1}{3}} & 0 & {\frac {1}{6}} & 0 & 0 & {\frac {1}{2}} & 0 & 0 & 0 & 0 & 0 & 0 & 0 & 0 & 0 \\
0 & 0 & {\frac {1}{2}} & {\frac {1}{2}} & 0 & 0 & 0 & 0 & 0 & 0 & 0 & 0 & 0 & 0 & 0 & 0 \\
{\frac {1}{6}} & 0 & 0 & 0 & 0 & {\frac {1}{2}} & 0 & 0 & {\frac {1}{3}} & 0 & 0 & 0 & 0 & 0 & 0 & 0 \\
0 & {\frac {1}{6}} & 0 & 0 & {\frac {1}{6}} & 0 & {\frac {1}{3}} & 0 & 0 & {\frac {1}{3}} & 0 & 0 & 0 & 0 & 0 & 0 \\
0 & 0 & {\frac {1}{6}} & 0 & 0 & {\frac {2}{9}} & {\frac {2}{9}} & {\frac {1}{6}} & 0 & 0 & {\frac {2}{9}} & 0 & 0 & 0 & 0 & 0 \\
0 & 0 & 0 & {\frac {1}{6}} & 0 & 0 & {\frac {1}{2}} & 0 & 0 & 0 & 0 & {\frac {1}{3}} & 0 & 0 & 0 & 0 \\
0 & 0 & 0 & 0 & {\frac {1}{3}} & 0 & 0 & 0 & 0 & {\frac {1}{2}} & 0 & 0 & {\frac {1}{6}} & 0 & 0 & 0 \\
0 & 0 & 0 & 0 & 0 & {\frac {2}{9}} & 0 & 0 & {\frac {1}{6}} & {\frac {2}{9}} & {\frac {2}{9}} & 0 & 0 & {\frac {1}{6}} & 0 & 0\\
0 & 0 & 0 & 0 & 0 & 0 & {\frac {1}{3}} & 0 & 0 & {\frac {1}{3}} & 0 & {\frac {1}{6}} & 0 & 0 & {\frac {1}{6}} & 0 \\
0 & 0 & 0 & 0 & 0 & 0 & 0 & {\frac {1}{3}} & 0 & 0 & {\frac {1}{2}} & 0 & 0 & 0 & 0 & {\frac {1}{6}} \\
0 & 0 & 0 & 0 & 0 & 0 & 0 & 0 & {\frac {1}{2}} & 0 & 0 & 0 & {\frac {1}{2}} & 0 & 0 & 0 \\
0 & 0 & 0 & 0 & 0 & 0 & 0 & 0 & 0 & {\frac {1}{2}} & 0 & 0 & {\frac {1}{6}} & 0 & {\frac {1}{3}} & 0 \\
0 & 0 & 0 & 0 & 0 & 0 & 0 & 0 & 0 & 0 & {\frac {1}{2}} & 0 & 0 & {\frac {1}{3}} & 0 & {\frac {1}{6}} \\
0 & 0 & 0 & 0 & 0 & 0 & 0 & 0 & 0 & 0 & 0 & {\frac {1}{2}} & 0 & 0 & {\frac {1}{2}} & 0
\end{array} \right)  \label{eq:M}
\end{equation}
You can in fact check that this is a Markov matrix, as all rows sum up to 1.  From the general theory of Markov chains [30] we know that at least one eigenvalue of this matrix has modulus 1 and that the moduli of all eigenvalues are no more than 1.  Moreover, the left eigenvector $v = \left[v_1, v_2, \ldots, v_{16} \right]$ corresponding to each eigenvalue equal to 1 solves the `balance equations', i.e.
$$\begin{array}{ccc}
M_{11} v_1 + M_{21} v_2 + \ldots + M_{16,1} v_{16} & = & v_1 \\
M_{12} v_1 + M_{22} v_2 + \ldots + M_{16,2} v_{16} & = & v_2 \\
\vdots & \vdots & \vdots \\
M_{1,16} v_1 + M_{2,16} v_2 + \ldots + M_{16,16} v_{16} & = & v_{16}
\end{array}$$
Thus, such an eigenvector represents the invariant measure of the Markov chain.  Specifically, for the Markov transition matrix $M$ shown in (\ref{eq:M}), this eigenvector is given by
\begin{eqnarray*}
& & \left[0.0159, 0.0466, 0.0648, 0.0368, 0.0481, 0.1013, 0.1182, 0.0445, \right. \\
& & \left. 0.0692, 0.1223, 0.1285, 0.0442, 0.0417, 0.0548, 0.0476, 0.0154 \right].
\end{eqnarray*}
Returning to the more readily interpretable representation in the $\{0,1,2,3\} \times \{0,1,2,3\}$ state space, the resulting invariant measure is given by
$$\left(\begin{array}{cccc}
0.0368 & 0.0445 & 0.0442 & 0.0154 \\
0.0648 & 0.1182 & 0.1285 & 0.0476 \\
0.0466 & 0.1013 & 0.1223 & 0.0548 \\
0.0159 & 0.0481 & 0.0692 & 0.0417
\end{array}\right)$$ 
where the horizontal axis corresponds to the number of sites that are $+1$, from $0$ to $3$, while the vertical axis corresponds to the number of arcs that are $+1$, from $0$ to $3$.  This means that, for example, the steady state probability of finding the system in the state $(1,2)$, i.e. with one site equal to $+1$ and two arcs equal to $+1$ is equal to $0.1182$.  The most likely outcomes of this Markov chain are the states $(2,1)$ and $(2,2)$, while the states $(0,0)$ and $(3,3)$ are the least likely ones.

\begin{figure}[htbp] 
   \centering
   \includegraphics[width=7in]{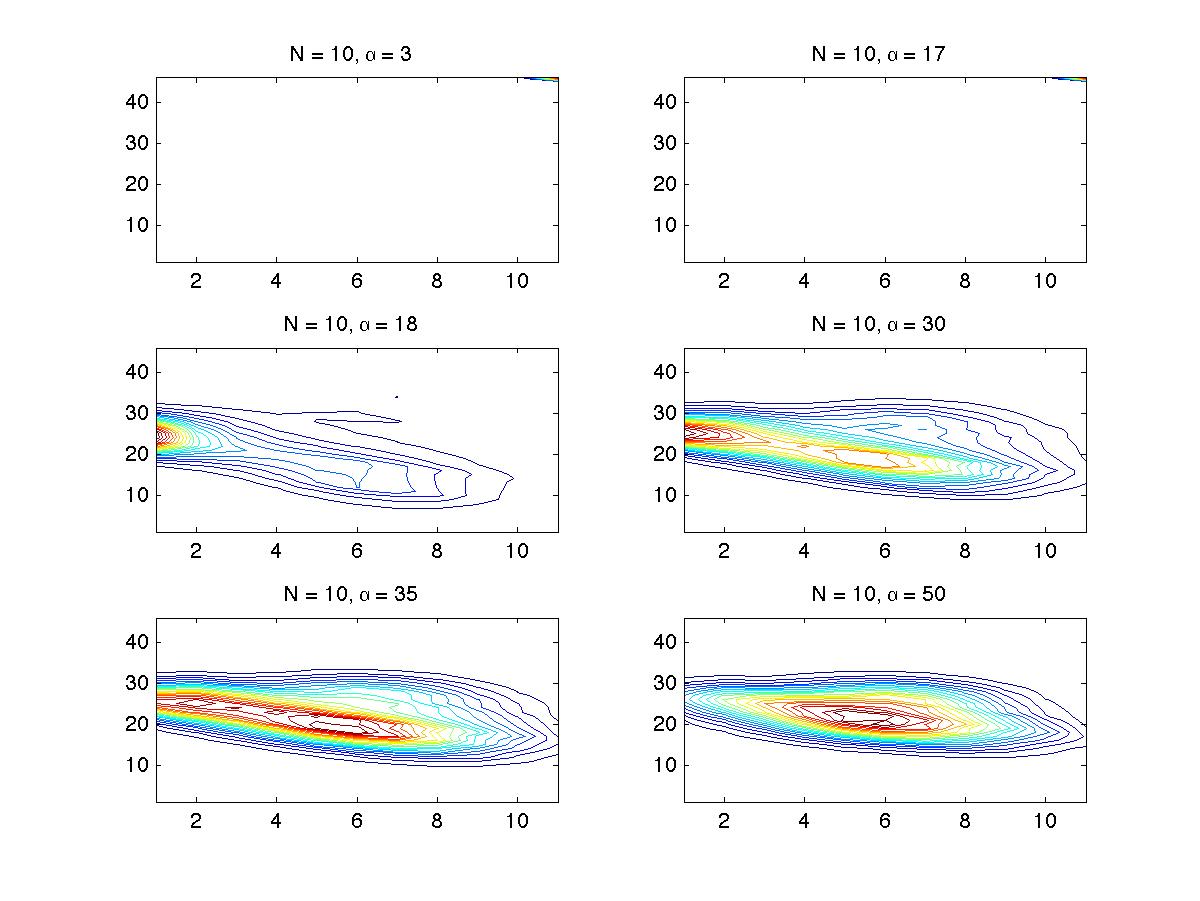} 
   \caption{Invariant measures for $N = 10$ and $\alpha = 3$ in the subcritical regime to $\alpha = 50$ in the supercritical regime.  The critical transition for $N = 10$ happens between $\alpha = 17$ and $\alpha = 18$.  As $\alpha$ increases in the supercritical regime, the invariant measure evolves from unimodal to bimodal and back to unimodal, with the mode moving from the left edge towards the middle of the state space.}
   \label{fig:invmeasure1}
\end{figure}

\begin{figure}[htbp] 
   \centering
   \includegraphics[width=4in]{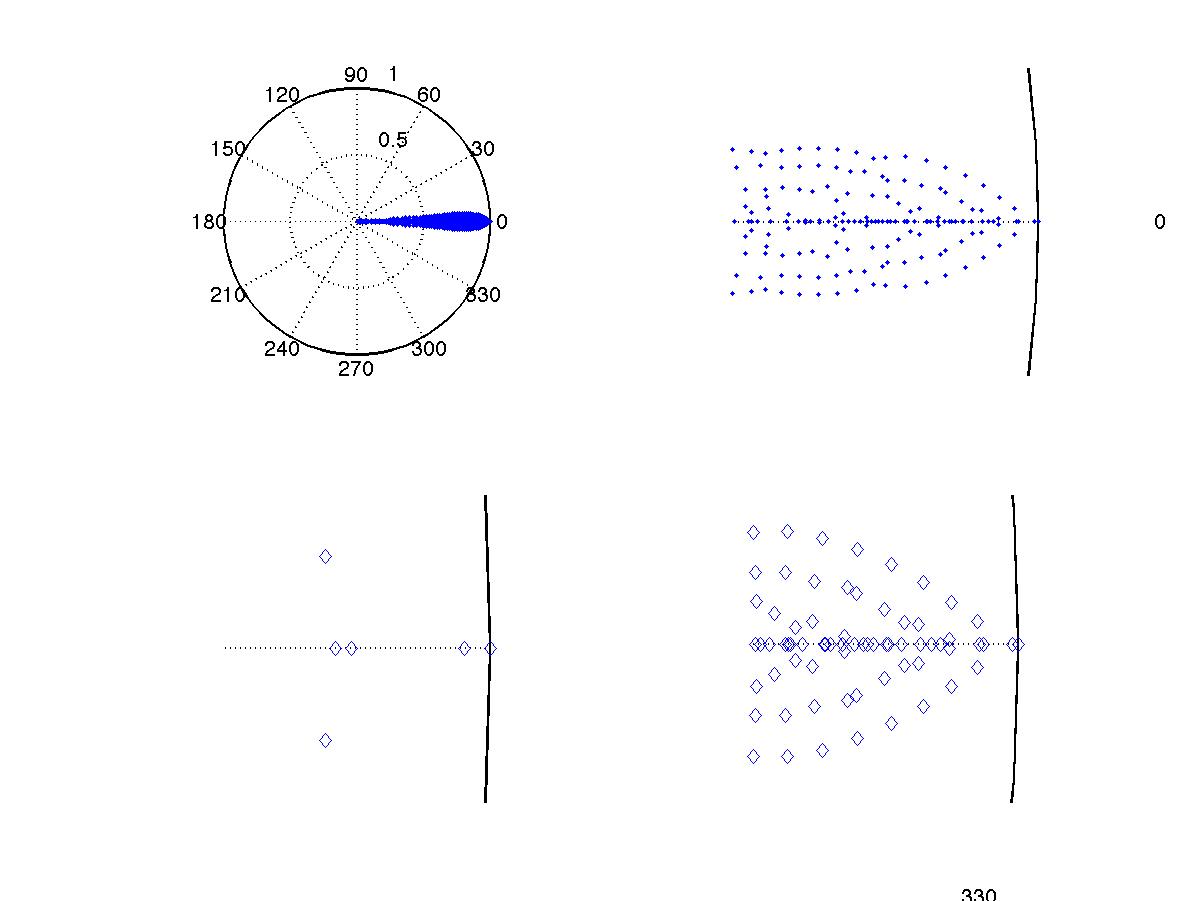} 
   \caption{The locations of the eigenvalues of the transition matrix for $N = 10$ and $\alpha  = 3$ in the unit disc.  The four panels progressively zoom in to the spectral gap in a clockwise manner.  The spectral gap is indeed very small, only equal to $1 - \lambda_2 \approx 1.3765 \times 10^{-6}$.}
   \label{fig:spectrum2}
\end{figure}

Appendix 3 describes the technical details involved in the appropriate convergence concepts for a Markov chain.  This discussion substantiates our use of linear algebra techniques to obtain information about the distributional path properties of $S^+$ and $A^+$.   We now proceed to investigate the changes to the resulting invariant measure as $\alpha$ and $N$ are allowed to vary.  Figure \ref{fig:invmeasure1} shows how the invariant measure evolves for $N = 10$ as $\alpha$ is allowed to increase.
As a first observation, note that the case $N = 10$ and $\alpha = 3$, whose deterministic approximation gave rise to the four separate attractors we described in section 4, possesses in fact a unique invariant measure, concentrated entirely at $\left(N,{\rm C}^N_2\right)$.  In fact, we can easily verify, using (\ref{eq:MCtransitions}), that when $N = 10$ and $\alpha = 3$, 
\begin{eqnarray*}
{\rm Pr} \left( \left. {\bf X}_{n+1} = (10,45) \right| {\bf X}_n = (10,45) \right) & = & {\rm Pr} \left( \left. {\bf X}_{n+1} = (10,45) \right| {\bf X}_n = (9,45) \right) \\
& = & {\rm Pr} \left( \left. {\bf X}_{n+1} = (10,45) \right| {\bf X}_n = (10,44) \right) = 1,
\end{eqnarray*}
thereby guaranteeing that $(10,45)$ is a trapping state for the full stochastic dynamics.  The only question that remains is one of ergodicity: might the corresponding Markov transition matrix possess a double eigenvalue equal to 1, in which case the two resulting eigenvectors will give rise to two distinct invariant measures, and the dynamics will be non-ergodic?
\begin{figure}[htbp] 
   \centering
   \includegraphics[width=4in]{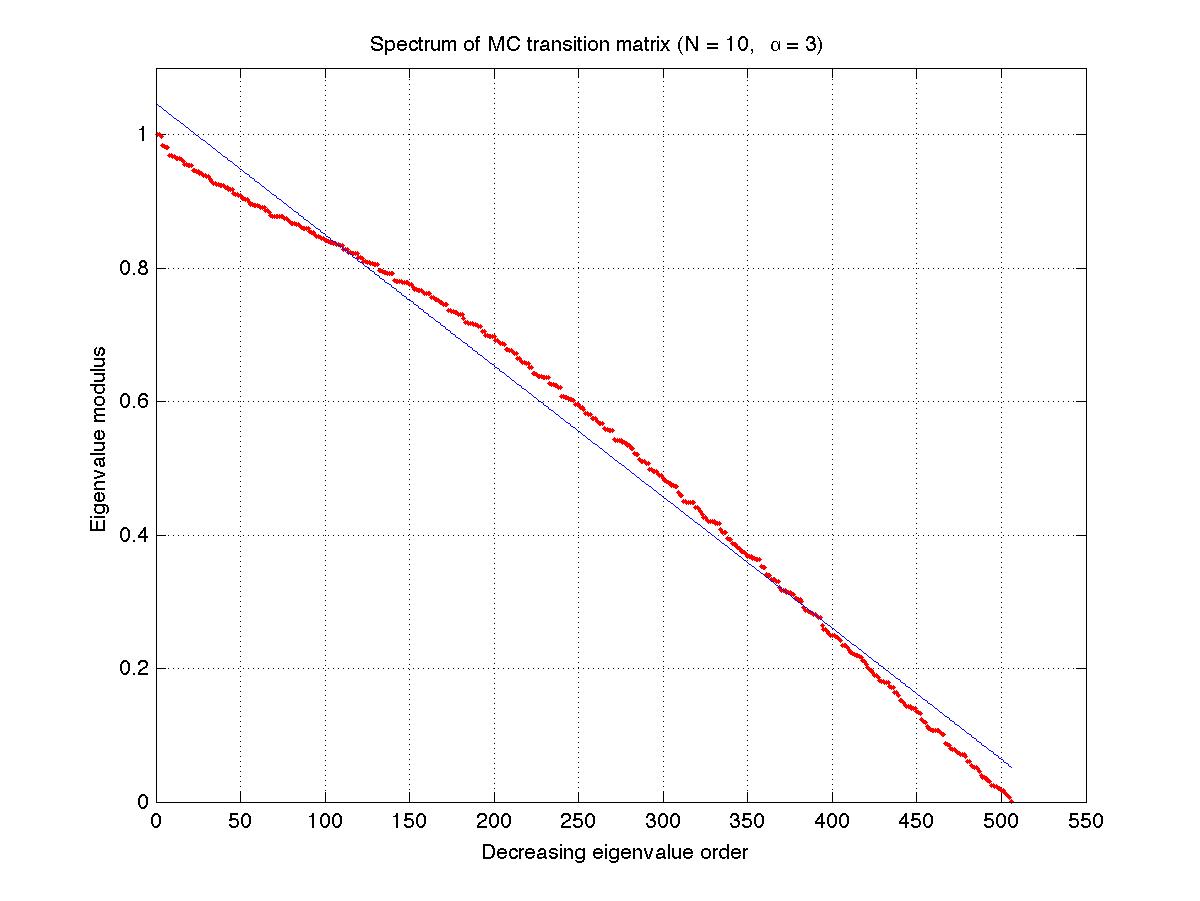} 
   \caption{The modulus of the eigenvalues of the transition matrix for $N = 10$ and $\alpha  = 3$ decay almost linearly, with $\lambda_k \approx -0.002k + 1.0466$.}
   \label{fig:spectrum1}
\end{figure}
Figure \ref{fig:spectrum1} shows the rate of decay of the eigenvalues for the transition matrix in this case.  Figure \ref{fig:spectrum2} shows the pattern generated by the location of the eigenvalues in the unit disc.  It turns out that the second highest eigenvalue is extremely close to $1$, at about $0.9999986235$!

\section{Transient Domain}
Thus, we know that the invariant measure is unique, despite the multiplicity of attractors for the deterministic cellular automaton approximation.  Why is that, and what can we infer from the nearness of the second eigenvalue to 1?  There are two justifications for the discrepancy between these two apparently contradictory asymptotic behaviors.  The first hinges on the observation, which we mentioned at the end of section 3, that the deterministic cellular automaton doesn't in fact provide a skeleton around which the Markov chain evolves.  Instead, these two dynamics deviate from one another systematically, in very instructive ways:
\begin{itemize}
\item Propagating the average is not the same as the average of the propagated measure.
\item Rounding the state to the nearest integer doesn't respect the sign of the expected increments.
\end{itemize}
The sources of these discrepancies are discussed in detail in Appendix 4.

Despite these important ways in which the underlying stochastic process deviates meaningfully from the deterministic cellular automaton approximation, in practice economic agents often use such approximation schemes. Specifically, stochastic dynamics are not uncommonly interpreted `quasi-deterministically', taking one step at a time, depending on the sign of signals on average, rather than propagate the ensemble of paths with their respective probabilities and subsequently take the average.

We saw the first justification for the deviation in asymptotic behaviors of the Markov chain and the deterministic cellular automaton.  According to it, we really didn't have much reason to expect that they would approximate one another in the first place!  The second justification for their deviation instead focuses on their hidden similarities, delving deeper into the interpretation of the spectrum of the transition matrix.  Specifically, as documented even in Figures \ref{fig:MC1} and \ref{fig:MC2}, the Markov chain with $N = 10$ and $\alpha  = 3$ doesn't converge to its invariant measure for a long time: in the first figure it took about $9,000$ steps, and in the second figure it didn't converge until after $10,000$ steps.  Moreover, in both cases the Markov chain spent most of the time oscillating in the vicinity of two states, around $(0,22)$ and $(5,8)$, occasionally switching between them.  We have already seen in Table 1 that these are in fact stable periodic attractors of the deterministic cellular automaton.  So the cellular automaton, despite being a lousy asymptotic approximation, is evidently quite accurate in describing the time evolution of particular paths of the stochastic process.  Why is this?

In fact, the reason for the transient, though long-lasting, relevance of the cellular automaton is the exceedingly slow converge rate of the Markov chain for $N = 10$ and $\alpha = 3$, which is due to its near-non-ergodicity.  As we already saw, the second highest eigenvalue is roughly within a millionth of $1$.  If it were equal to $1$, then the Markov chain would lose ergodicity, and the asymptotic multiplicity of attractors exhibited by the cellular automaton would dominate the asymptotic behavior of the stochastic system as well.  As it is, this is averted, but only by the slightest of margins.  In fact, the convergence rate of an ergodic Markov chain to its unique invariant measure is controlled exponentially by the second highest eigenvalue, $\lambda_2$.  More specifically [30], for any function $F: {\mathcal S} \times {\mathcal A} \rightarrow {\Bbb R}$, the variance of the expectation of $F$ under the measure $\pi_n$ falls with $n$ as $\lambda_2^n$.  In the case of $N = 10$ and $\alpha = 3$, this remains high for a very long time.

For comparison, $\lambda_2 = 0.9942$ when $N = 10$ and $\alpha  = 18$.  This would take $120$ steps to half the variation distance\footnote{i.e. the distance 
$$\delta\left(\mu^{(n)},\mu^{(\infty)} \right) = \sup_k \left| \mu_k^{(n)} - \mu_k^{(\infty)} \right|$$} from the invariant measure, while the $N = 10$, $\alpha = 3$ case we considered above would take $503,558$ steps to reach the same approach to the invariant measure!  As a result of this effect, we ever only get to experience the transient regime of the Markov chain for $N = 10$ and $\alpha = 3$.    
\begin{figure}[htbp] 
   \centering
   \includegraphics[width=3.5in]{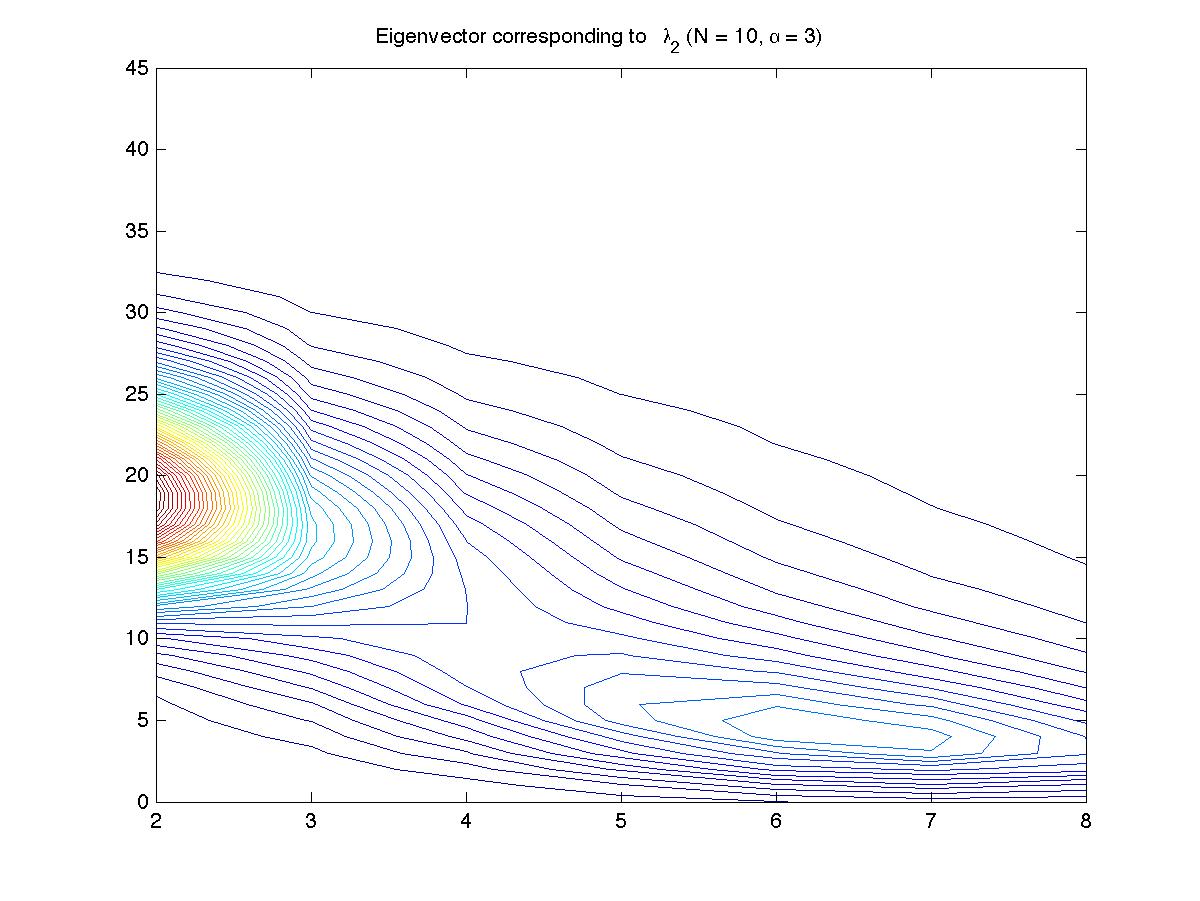} 
   \caption{Contour plot for the eigenvector corresponding to $\lambda_2$ when $N = 10$ and $\alpha = 3$.  The axes have been cropped because the eigenvector is dominated by a positive peak in the neighborhood of $(0,22)$ and a negative peak at $(10,45)$.  This figure shows that there is a ridge connecting the vicinity of $(0,22)$ and that of $(5,8)$.}
   \label{fig:spectrum3}
\end{figure}
\begin{figure}[htbp] 
   \centering
   \includegraphics[width=6.5in]{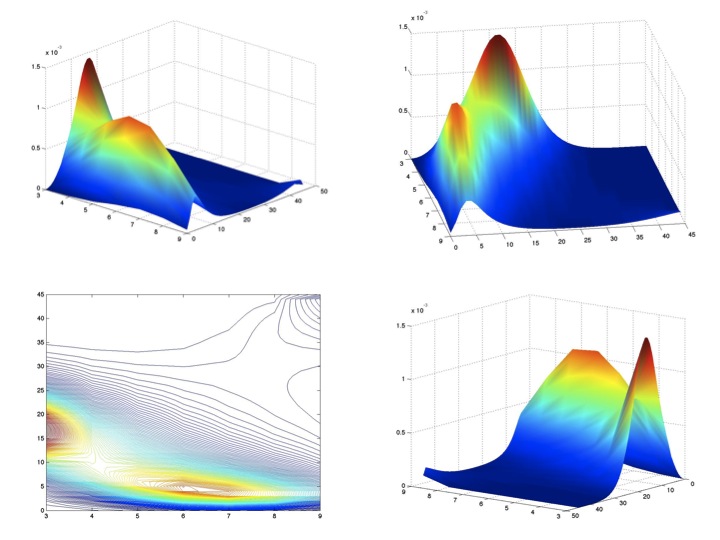} 
   \caption{3D visualization of the eigenvector corresponding to $\lambda_2$ for $N = 10$ and $\alpha = 3$.  At this cropping, the two peaks that dominate the transient behavior are clearly recognizable, as is the ridge the connects them, along which the system occasionally transitions.  At the corner near $(10,45)$, which constitutes the peak of the invariant measure and has been cropped out, we are beginning to discern a shallow slope, which represents the last barrier before the trapping state at $(10,45)$ is reached.  The lower left panel illustrates the contour lines of the 3D figures, clearly depicting the paths along which the process will eventually become channeled towards the invariant measure.}
   \label{fig:spectrum4}
\end{figure}
Figures \ref{fig:spectrum3} and \ref{fig:spectrum4} illustrate the information gleaned from the eigenvector associated with the second highest eigenvalue, $\lambda_2$.  These visualization tools substantiate our view of the Markov chain `mixing' the attractors of the deterministic cellular automaton approximation.  Rather than irrelevant, the cellular automaton, and the contingent submartingale method that motivated it, cast useful light on the transient behavior of the Markov chain, which in this instance is arguably more relevant than its invariant measure, which won't be reached for an inordinate amount of time.   

\section{Phenomenology of Distributional Properties}
\begin{figure}[htbp] 
   \centering
   \includegraphics[width=7in]{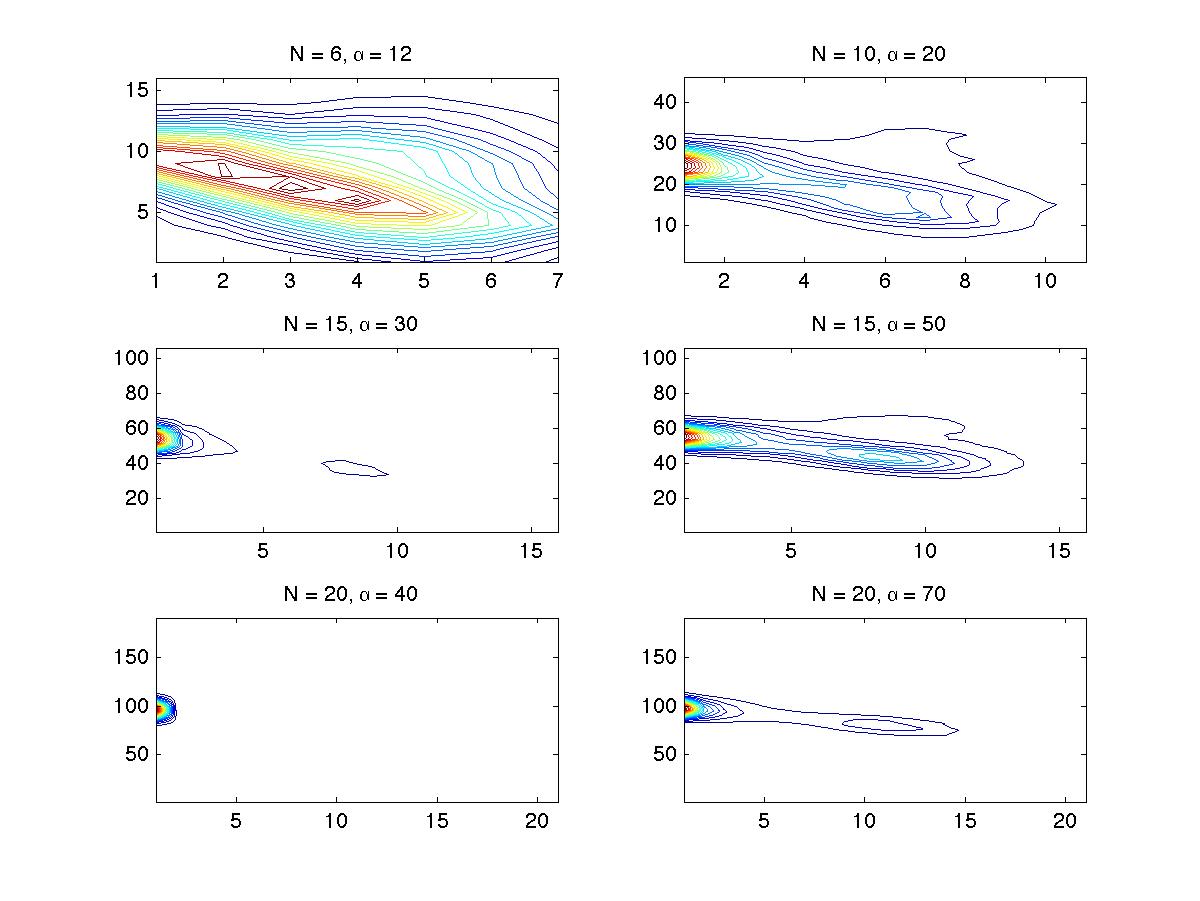} 
   \caption{Invariant measures for $N = 6, 10, 15, 20$ and $\alpha$ in the supercritical regime.  The second critical transition appears to occur for $\alpha > 2N$ as $N$ increases.  But when it occurs, it appears to follow the same pattern of evolution.}
   \label{fig:invmeasure2}
\end{figure}

Armed with the techniques and insights form the previous two sections, we are now in a position to describe the qualitative changes that the invariant measure experiences as the two parameters $N$ and $\alpha$ are varied.  
\begin{enumerate}
\item[1)] For every $N$ there exists a critical value of the coupling constant $\alpha^\ast$ that separates the sub-critical from the super-critical regime of the Markov chain.  This is analogous to the critical phase transition described in [32].  In this instance, the subcritical regime entails a unique invariant measure concentrated at the trapping state $\left( N,{\rm C}^N_2 \right)$.  The super-critical regime on the other hand entails a unique distribution with broad support.  More specifically, we conjecture that for all $N$ there exists a value $\alpha^\ast$ such that for all $\alpha < \alpha^\ast$, $\pi_\infty \left(N,{\rm C}^N_2 \right) = 1$ and for all $\alpha \geq \alpha^\ast$, 
$${\bf E}^{\pi_\infty} \left[ S^+ A^+ \right] \ne {\bf E}^{\pi_\infty} \left[ S^+ \right] {\bf E}^{\pi_\infty} \left[ A^+ \right].$$
\item[2)] As corroborated by Figure \ref{fig:invmeasure1}, the critical value $\alpha^\ast$ is lower than $2N$.  In fact we conjecture that $\alpha^\ast = 2N-2$. 
\item[3)] As $\alpha$ increases past $\alpha^\ast$, the invariant measure is initially bimodal, a dominant one in the vicinity of $\left( 0,{\frac {{\rm C}^N_2}{2}} \right)$ and a substantially shallower one around $\left( {\frac {N}{2}},{\frac {{\rm C}^N_2}{6}} \right)$.  We conjecture that 
$$\arg \min_{{\mathcal S} \times {\mathcal A}} \pi_\infty \left(i,j;\alpha^\ast \right) \in \{0,1\} \times \left\{ \left\lfloor {\frac {{\rm C}^N_2}{2}} \right\rfloor, \left\lceil {\frac {{\rm C}^N_2}{2}} \right\rceil \right\}.$$
\item[4)] Progressively, as $\alpha$ is increased, mass is transferred from the former mode to the latter one, and the latter mode slowly slides up towards $\left( {\frac {N}{2}},{\frac {{\rm C}^N_2}{2}} \right)$.  Asymptotically as $\alpha \rightarrow \infty$ the invariant measure remains broadly supported, and becomes symmetric around a unique peak in the vicinity of $\left( {\frac {N}{2}},{\frac {{\rm C}^N_2}{2}} \right)$.  Specifically, we conjecture that there exists a $\sigma$, such that
$$\lim_{\alpha \rightarrow \infty} \pi_\infty (i,j;\alpha) = {\frac {1}{2 \pi \sigma^2}} \exp \left\{-{\frac {4(2i-N)^2 + \left(4j-N^2+N \right)^2}{32 \sigma^2}} \right\}.$$
\item[5)] As $N$ increases, the eventual mixing in the super-critical regime requires larger increases in $\alpha$.  This effect is illustrated in Figure \ref{fig:invmeasure2}.  We conjecture that for every $N$ there exists a second critical value of the coupling constant, $\alpha^{**} (N)$, such that for all $\alpha \geq \alpha^{**}$, $\pi_\infty$ becomes unimodal once more.  Moreover, we conjecture that 
$$\lim_{N \rightarrow \infty} {\frac {\log \left(\alpha^{**} (N) \right)}{\log N}} > 1.$$
\item[6)] For every $N$, the spectral gap is an increasing function of $\alpha$.  Specifically we conjecture that it is a power law, with two branches, separated by $\alpha^\ast$, i.e. that there exist positive constants $C_1$ and $C_2$ such that for any $0<\alpha_1<\alpha_2<\alpha^\ast<\alpha_3<\alpha_4$,
$${\frac {\log \left({\frac {1-\lambda_2(\alpha_2)}{1-\lambda_2(\alpha_1)}} \right)}{\log \left( {\frac {\alpha_2}{\alpha_1}} \right)}} = C_2 > C_1 = {\frac {\log \left({\frac {1-\lambda_2(\alpha_4)}{1-\lambda_2(\alpha_3)}} \right)}{\log \left( {\frac {\alpha_4}{\alpha_3}} \right)}}.$$ 
\item[7)] Returning to the sub-critical regime we saw that the behavior of the Markov chain is dominated by the transient regime because of extremely high $\lambda_2$.  Moreover we saw that the attractors of the deterministic cellular automaton approximation coincided with the structure of the eigenvector corresponding to $\lambda_2$, guiding the Markov chain to occasionally transition among the stable attractors of the cellular automaton.  The inexorable convergence to the trapping state doesn't materialize until times arguably beyond practical relevance.  We conjecture that this behavior is characteristic of the sub-critical regime. 
\end{enumerate}

Finally, we proceed to offer a speculative interpretation of the distributional properties of this Markov chain  that represents the interactions of hierarchical economic agents.  In order to draw economically relevant conclusions, we introduce the concepts of {\it connectivity}, {\it conservatism} and {\it conformity}, along with the opposite respective poles, {\it fragmentation}, {\it liberalism} and {\it individualism}.  In our model, we interpret the site spins as acceptance ($+1$) or rejection ($-1$) of a new idea (along the conservatism to liberalism dimension), and the arc spins as acceptance ($+1$) or rejection ($-1$) of a social connection that enhances cohesion.  Finally, we interpret the coupling constant $\alpha$, which serves to modulate the tradeoff between the pursuit of local majority and that of global minority, as a measure of conformity: low $\alpha$ represents societies where people prioritize conforming to group norms, while high values of $\alpha$ represent societies where people prioritize seizing opportunities for speculative contrarian actions.

With this interpretative framework in mind, we can draw the following three conclusions, inspired by the qualitative properties of the distributions that arise from our interacting hierarchical agents:
\begin{enumerate}
\item[(I)] Societies with a high desire for uniformity are driven towards broad, multilateral connectivity and uniform liberalism.  This represents the trapping state equilibrium in the sub-critical regime of our model.
\item[(II)] As the relevance of individualism increases in a society, we see the highest structural diversity in outcomes, dominated by uniform conservatism and centralized connectivity.  This represents intermediate values of $\alpha$, just past the phase transition, and the resulting multimodal distributions.
\item[(III)] Societies with a high degree of individualistic opportunism, become tolerant of diversity of opinions and fragmented into local islands of social cohesion.  This represents the eventually unimodal, broadly supported distributions we obtain asymptotically as $\alpha$ increases further.
\end{enumerate}

This scheme for interpreting the qualitative conclusions of our study is but a caricature.  Nevertheless, they bear some striking similarities with socioeconomic systems in various historical and contemporary settings.  It may be the case that further elaboration of these parallels between our model and socioeconomic concepts will substantiate these tentative suggestions.  Research in this direction appears warranted.

\section{Conclusions and Next Steps}
In summary, we defined a stochastic process on meta-individual agents that we represent as spin configurations of sites and arcs in a network that evolves as a result of their interactions.  As such, this work can be seen as extending earlier work on spin market models and opinion formation models to the case of a network that evolves simultaneously and reciprocally with the agents' actions.

We were able to obtain analytic control of the invariant measure for the resulting stochastic process, and compare it with a deterministic cellular automaton that was designed to approximate it.  Our analysis shed light on questions of ergodicity and persistent transient effects.  Based on the novel analytic methodology we describe, a series of conjectures is offered, to help guide future work in this field.

Beyond the specific technical questions that are covered by these conjectures, larger next steps include the extension of the framework proposed here to more layers in the agents' hierarchy, including objects beyond arcs, for instance `faces' or triples of sites.  Finally, another inviting direction for future work involves delineating the speculative socioeconomic interpretations offered at the end of the previous section.  It would be worthwhile to attempt a historical/empirical investigation of the proposed parallels.  For now, we must settle for the many more questions than answers that our study brought to light so far.

\section*{Acknowledgement}
The author would like to thank two anonymous referees for recommending several useful references and for suggestions that helped streamline the clarity of the presentation.


\section{References}
\begin{enumerate}
\item[[1]] Alfarano S., Lux T.  A noise trader model as a generator of apparent financial power laws and long memory. {\it Macroeconomic Dynamics}, 11:80-101, 2007.

\item[[2]] Alfarano S., Wagner F., Lux T.  Estimation of agent-based models: the case of an asymmetric herding model. {\it Computational Economics}, 26:19-49, 2005.

\item[[3]] Aoki M.  {\it Modeling aggregate behavior and fluctuations in economics: stochastic views of interacting agents}.  Cambridge University Press, Cambridge, 2004.

\item[[4]] Aoki M.  Thermodynamic limits of macroeconomic or financial models: one-parameter and two-parameter Poisson-Dirichlet models.  {\it Journal of Economic Dynamics and Control}, 32:66-84, 2008.

\item[[5]] Barnsley M.F., Demko S.G., Elton J.H., Geronimo J.S., Invariant measures for markov processes arising from iterated function systems with place-dependent probabilities.  {\it Annales de l'Institut Henri Poincar\'e, Probabilit\'es et Statistiques}, 24:367-394, 1988.

\item[[6]] Benediktsdottir S.  An empirical analysis of specialist trading behavior at the New York Stock Exchange.  {\it International Finance Discussion Papers, Board of Governors of the Federal Reserve System}, (876), 2006.

\item[[7]] Bouchaud J.-P., M\'ezard M., Potters  M., Statistical properties of stock order books and price impact.  {\it Quantitative Finance}, 2:251-256, 2002.

\item[[8]] Brock W.A., Hommes C.H., Heterogeneous beliefs and routes to chaos in a simple asset pricing model.  {\it Journal of Economic Dynamics and Control}, 22:1235-1274, 1998.

\item[[9]] Chiarella C., Iori G., Perell\'o J., The impact of heterogeneous trading rules on the limit order book and order flows.  {\it Journal of Economic Dynamics and Control}, 33:525-553, 2009.

\item[[10]] Cioffi-Revilla C., Computational social science.  {\it Wiley Interdisciplinary Reviews: Computational Statistics}, 2:259-271, 2010.

\item[[11]] Cirillo P., Redig F., Ruszel W.,  Duality and stationary distributions of wealth distribution models.  {\it Journal of Physics A: Mathematical and Theoretical}, 47:085203, 2014.

\item[[12]] Durrett R.  Ten lectures in probability theory.  In {\it Lectures in Probability Theory (Saint-Flour, 1993)}, volume 1608 of {\it Lecture Nots in Mathematics}, pages 97-201.  Springer, 1995.

\item[[13]] Farmer J.D., Gillemot L., Lillo F., Mike S., Sen A., What really causes large price changes?.  {\it Quantitative Finance}, 4:383-397, 2004.

\item[[14]] Feller W.  {\it An introduction to probability theory and applications, vol II}.  Wiley, New York, 1971.

\item[[15]] F\"ollmer H., Horst U., Kirman A.,  Equilibria in financial markets with heterogeneous agents: a probabilistic perspective.  {\it Journal of Mathematical Economics}, 41:123-155, 2005.

\item[[16]] Gabaix X., Gopalakrishnan P., Plerou V., and Stanley H.E.  A theory of power-law distributions in financial fluctuations.  {\it Nature}, 423:267-270, 2003.

\item[[17]] Ghoulmie F., Cont R., Nadal J.-P.,  Heterogeneity and feedback in an agent-based market model.  {\it Journal of Physics: Condensed Matter}, 17:S1259-S1268, 2005.

\item[[18]] Gillemot L., Farmer J.D., Lillo F.,  There's more to volatility than volume.  {\it Quantitative Finance}, 6:371-384, 2006.

\item[[19]] Horst U.  Dynamic systems of social interactions.  {\it Journal of Economic Behavior and Organization}, 73: 158-170, 2010.

\item[[20]] Horst U., Rothe C.,  Queuing, social interactions and the microstructure of financial markets.  {\it Macroeconomic Dynamics}, 12:211-233, 2008.

\item[[21]] Johnson N.,Jefferies P., and Hui P.  {\it Financial market complexity}.  Oxford University Press, 2003.

\item[[22]] Kaizoji T., Bornholdt S., and Fujiwara Y.  Dynamics of price and trading volume in a spin market model of stock markets with heterogeneous agents.  {\it Physica A}, 316:441-452, 2002.

\item[[23]] Kirman A.  Whom or what does the representative individual represent?. {\it Journal of Economic Perspectives}, 6:117-136, 1992.

\item[[24]] Lo A.,  Long-term memory in stock market prices.  {\it Econometrica}, 59:1279-1313, 1991.

\item[[25]] Mandelbrot B.,  Statistical methodology for non periodic cycles: from the covariance to r/s analysis.  {\it Annals of Economic and Social Measurement}, 1(3):259-290, 1972.

\item[[26]] Potters M., Bouchaud J.P.,  More statistical properties of stock order books: empirical results and models.  {\it Physica A}, 324:133-140, 2003.

\item[[27]] Singer I.M. and Thorpe J.A.  {\it Lecture notes in elementary topology and geometry}.  Springer, 1976.

\item[[28]] Steinsaltz D., Locally contractive iterated function systems.  {\it Annals of Probability}, 27:1952-1979, 1999.

\item[[29]] Stroock D.  {\it Markov Processes from K. It\^o's perspective}. Princeton University Press, Princeton, 2003.

\item[[30]] Stroock D.  {\it An introduction to Markov processes}. Springer, Berlin, 2005.

\item[[31]] Tedeschi G., Iori G., Gallegati M., The role of communication and imitation in limit order markets.  {\it European Physical Journal B}, 1:489-497, 2009.

\item[[32]] Theodosopoulos T.  Uncertainty relations in models of market microstructure.  {\it Physica A}, 355:209-216, 2005.

\item[[33]] Theodosopoulos T. and Yuen M.  Properties of the wealth process in a market microstructure model,  {\it Physica A}, 378: 443-452, 2007.

\item[[34]] Wolfram S.  {\it A new kind of science}.  Wolfram Media, Champaign, 2002.

\end{enumerate}

\section*{Appendix 1: Neighborhood Structures}
This Appendix provides more technical details about the construction of the random neighborhoods, ${\mathcal N}_\vartheta (y,T_n)$ and ${\mathcal N}_\eta \left(x,y,T_n \right)$, for a site $y \in {\mathcal S}$ and arc $(x,y) \in {\mathcal A}$.

In particular ${\mathcal N}_\vartheta (y,T_n)$ is treated as a random subset of ${\mathcal S} \setminus \{y\}$ created as follows: 
\begin{quotation}
\noindent Represent the members of ${\mathcal A}$ as ${\rm C}_2^N = \left( \begin{array}{c} N \\ 2 \end{array} \right)$ balls in an urn, with $A^+ \left(T_n \right)$ of them WHITE and the remaining BLACK.  Draw $N-1$ balls from this urn, representing the members of ${\mathcal S} \setminus \{y\}$.  The WHITE balls that we draw represent the members of the neighborhood ${\mathcal N}_\vartheta (y,T_n)$.
\end{quotation}
Thus, ${\mathcal N}_\vartheta (y,T_n)$ is random subset of ${\mathcal S} \setminus \{y\}$ chosen uniformly among all those with cardinality ${\bf V}(y,T_n)$, a hypergeometric random variable, with $N-1$ draws without replacement from a population of ${\rm C}_2^N$ with $A^+(T_n)$ successes, and $\{ \{{\bf V} \}_{n=0}^\infty \}_{y \in {|mathcal S}}$ is an independent family of random variables.

In the case of the arc configurations, ${\mathcal N}_\eta \left(x,y,T_n \right)$ is treated as a random subset of ${\mathcal A} \setminus \{(x,y)\}$ created as follows:
\begin{quotation}
\noindent Represent the members of ${\mathcal S}$ as $N$ balls in an urn, with $S^+\left( T_n \right)$ of them RED and the remaining BLUE.  Draw two balls from this urn, representing the two endpoints of the base arc $(x,y)$.  The number of RED balls that we draw represent the endpoints of the base arc that accept neighbor connections.  All $N-2$ arcs coming out of each of those chosen endpoints (i.e. the arcs to all but the two points involved in the base arc) taken together represent the neighborhood ${\mathcal N}_\eta \left(x,y,T_n \right)$.
\end{quotation}
Thus, ${\mathcal N}_\eta \left(x,y,T_n \right)$ is a random member of 
$$\left\{ \emptyset, {_xA^y} \cup {^yA_x}, {_yA^x} \cup {^xA_y}, {_xA^y} \cup {^yA_x} \cup {_yA^x} \cup {^xA_y} \right\}.$$ 
With probability $\left( 1- {\frac {S^+(T_n)}{N}} \right) \left(1 - {\frac {S^+(T_n)}{N-1}} \right)$, both endpoints of the base arc have spin $-1$ and therefore the base arc has no neighbors, making ${\mathcal N}_\eta \left(x,y,T_n \right) = \emptyset$.  On the other hand, with probability ${\frac {2S^+(T_n) \left(N - S^+(T_n) \right) }{N(N-1)}}$ the base arc has one endpoint with spin $+1$ and the other with spin $-1$, making $\left|{\mathcal N}_\eta \left(x,y,T_n \right)\right| = N-2$.  Finally, with probability ${\frac {S^+(T_n) \left(S^+(T_n) -1 \right)}{N(N-1)}}$ both endpoints of the base arc have spins equal to $+1$, and therefore $\left|{\mathcal N}_\eta \left(x,y,T_n \right)\right| = 2N-4$.  Let ${\bf Z} \left(x,y,T_n \right)$ be a hypergeometric random variable, with $2$ draws without replacement from a population of $N$ with $S^+(T_n)$ successes, and $\{ \{{\bf Z} \}_{n=0}^\infty \}_{(x,y) \in {\mathcal A}}$ an independent family of random variables, and independent of the earlier family $\{ \{{\bf V} \}_{n=0}^\infty \}_{y \in {\mathcal S}}$.  The cardinality of ${\mathcal N}_\eta \left(x,y,T_n \right)$ is given by ${\bf W} = {\bf Z} (N-1)$.
\\
\\
\section*{Appendix 2: Hypergeometric Conditional Tails}
This Appendix describe the technical details of the derivation of the transition probabilities, $P_{++} (i,j)$, $P_{--} (i,j)$, $Q_{++} (i,j)$ and $Q_{--} (i,j)$.  In particular, our goal here is to show the origin of the conditionally hypergeometric random variables, whose partial sums are involved in the computation of the transition probabilities, as shown in Section 3.

We begin with $P_{++}(i,j)$, the probability that, having chosen a $+1$ site, it remains $+1$ after the update, assuming the system is in a state with $S^+ = i$ and $A^+ = j$.  This probability has two components.  First, we need to compute the probability that the randomly chosen $+1$ site has $\ell$ neighbors, and then the chance that $k$ of them are positive, in order to compute the `local majority' contribution to the interaction potential.  Using the hypergeometric distribution, we can see that the former is given by 
$$\left(\begin{array}{c} {\rm C}^N_2 \\ N-1 \end{array} \right)^{-1} \left(\begin{array}{c} j \\ \ell \end{array} \right) \left(\begin{array}{c} {\rm C}^N_2 -j \\ N-\ell-1 \end{array} \right),$$
while the latter is given by
$$\left(\begin{array}{c} N-1 \\ \ell \end{array} \right)^{-1} \left(\begin{array}{c} i-1 \\ k \end{array} \right) \left(\begin{array}{c} N-i \\ \ell-k \end{array} \right).$$

The only thing that remains is to figure out the limits of the desirable range of values for $\ell$ and $k$.  Our goal is to guarantee that  $h_{\eta \oplus \vartheta} (x)>0$, where $x \in Y$ is the chosen site, and since we've assumed that the site is positive, this amounts to enforcing the inequality
$$\sum_{y \in {\mathcal N}_\vartheta (x)} \eta(y) > {\frac {\alpha}{2}} \left|{\frac {4(i+j)}{N(N+1)}} - 1\right|.$$
Assuming that there are $k$ positive out of $\ell$ total neighbors for $x$ (given the state $\eta \oplus \vartheta$), we see that $\sum_{y \in {\mathcal N}_\vartheta (x)} \eta(y) = 2k-\ell$.  Thus, given $\ell$, we need $k$ to satisfy $k>{\frac {1}{2}} \left( \ell+ {\frac {\alpha}{2}} \left|{\frac {4(i+j)}{N(N+1)}} -1 \right| \right)$. 

We now turn our attention to the probability $Q_{++} (i,j)$ that, having chosen a $+1$ arc, it remains $+1$ after the update, assuming the system is in a state with $S^+ = i$ and $A^+ = j$.  As we saw in Appendix 1, the number of neighbors of an arc can be either $0$, $N-2$ or $2(N-2)$, depending on whether none, one or both of its endpoints are $+1$ sites respectively.  Clearly, the first case, in which both of the chosen arc's endpoints are $-1$ sites, doesn't contribute anything to this probability, because in that case the arc has no neighbors and therefore its interaction potential has no `local majority' component.  The remaining `global imbalance' term in the interaction potential always has the opposite sign from the current sign of the site or arc under consideration, so it always makes the chosen site or arc flip its sign.  

Focusing on the second case, in which the chosen arc has only one $+1$ site as an endpoint, we see that this occurs with probability ${\frac {2i(N-i)}{N(N-1)}}$.  The only question that remains is how many of the $N-2$ resulting neighbors to the chosen arc are themselves positive.  The probability that $k$ out of the $N-2$ arc neighbors are positive is given by the hypergeometric distribution as
$$\left(\begin{array}{c} {\rm C}^N_2 -1 \\ N-2 \end{array} \right)^{-1} \left(\begin{array}{c} j-1 \\ k \end{array} \right) \left(\begin{array}{c} {\rm C}^N_2 -j \\ N-k-2 \end{array} \right).$$

Similarly, in the final case, in which both endpoints of the chosen arc are $+1$ sites, and therefore the chosen arc has $2(N-2)$ neighbors, occurs with probability ${\frac {i(i-1)}{N(N-1)}}$.  Finally, the probability that $k$ out of these $2N-4$ arc neighbors are positive is computed again using the hypergeometric distribution as
$$\left(\begin{array}{c} {\rm C}^N_2 -1 \\ 2N-4 \end{array} \right)^{-1} \left(\begin{array}{c} j-1 \\ k \end{array} \right) \left(\begin{array}{c} {\rm C}^N_2 -j \\ 2N-k-4 \end{array} \right).$$

As before, our goal this time is that $g_{\eta \oplus \vartheta} (x,y) >0$, where $(x,y) \in {\mathcal A}$ is the chosen arc.  Since we've assumed that this chosen arc is positive, this condition is equivalent to the following inequality:
$$\sum_{(u,v) \in {\mathcal N}_\eta \left(x,y\right)} \vartheta(u,v) > {\frac {\alpha}{2}} \left|{\frac {4(i+j)}{N(N+1)}} - 1\right|.$$  
In the case of $N-2$ neighbors, $k$ of which are positive, this becomes 
$$k>{\frac {1}{2}} \left(N-2+ {\frac {\alpha}{2}} \left|{\frac {4(i+j)}{N(N+1)}} -1 \right| \right).$$  
On the other hand, when the chosen arc has $2N-4$ neighbors, $k$ of which are positive, the inequality becomes
$$k>{\frac {1}{2}} \left(2N-4+ {\frac {\alpha}{2}} \left|{\frac {4(i+j)}{N(N+1)}} -1 \right| \right).$$ 
\\
\\
\section*{Appendix 3: Weak Convergence Formalism for Markov Chains}
This Appendix provides the detailed formalism behind the convergence of the empirical measures of a Markov chain and the questions of ergodicity that ensue.   At this point we will introduce some notation to help clarify the nature of the relevant limiting behavior.  Let 
$$\pi_n (i,j;\alpha,\beta) = {\rm Pr} \left({\bf X}_n = (i,j) \right)$$
at inverse temperature $\beta$ and coupling constant $\alpha$.  Then
$$\pi_n (i,j;\alpha) = \lim_{\beta \rightarrow \infty} \pi_n (i,j;\alpha,\beta)$$
and
$$\pi_\infty (i,j;\alpha) = \lim_{n \rightarrow \infty} \pi_n (i,j;\alpha).$$
Let ${\bf M} \in {\Bbb R}^{(N+1)\left(N^2 - N + 2 \right)}$ be the Markov transition matrix and 
$$\mu \in {\mathcal M}_1 \left(\left\{1,2, \ldots, (N+1) \left({\rm C}^N_2 + 1\right) \right\} \right)$$ 
be a probability measure, represented as a row (i.e. dual, or co-) vector with $(N+1) \left({\rm C}^N_2 + 1\right)$ dimensions.  In particular, we can represent any function 
$$h: \left\{1,2, \ldots, (N+1) \left({\rm C}^N_2 + 1\right) \right\} \rightarrow {\Bbb R}$$ 
as a vector ${\bf u}$ with $(N+1) \left({\rm C}^N_2 + 1\right)$ dimensions, and then the expected value of this function with respect to the probability measure $\mu$ is given by the inner product\footnote{This description is inspired by Stroock's account of Doeblin's theory of Markov chains in [30].} 
$$\mu {\bf u} = \sum_{\ell = 1}^{(N+1) \left({\rm C}^N_2 + 1\right)} \mu_\ell {\bf u}_\ell.$$

The Markov chain entails the propagation of the row vectors $\mu$ by the transition matrix ${\bf M}$.  Thus, starting with some initial co-vector $\mu^{(0)}$, the Markov chain can be described as the sequence of co-vectors $\mu^{(n+1)} = \mu^{(n)} {\bf M}$.  However, this propagation is difficult to interpret in terms of our stochastic process ${\bf X}_n$, because the one-dimensional row vectors of $(N+1) \left({\rm C}^N_2 + 1\right)$ dimensions are actually representing two-dimensional distributions on a $(N+1) \times \left({\rm C}^N_2 + 1\right)$ grid.  The transformation between the 1D co-vector representation and the more readily interpretable representation as measures on this $(N+1) \times \left({\rm C}^N_2 + 1\right)$ grid is effected by the functions $I$ and $J$ introduced in (\ref{eq:Ifunction}) and (\ref{eq:Jfunction}) respectively.  Thus, formally, 
$$\mu^{(n)} = \pi_n \circ (I \times J),$$
i.e.
$$\mu^{(n)}_k = \pi_n \left(I(k),J(k);\alpha \right),$$
where $\circ$ denotes function composition.  Inverting this we obtain
$$\pi_n (i,j;\alpha) = \mu^{(n)} \left(j\left[ {\rm C}^N_2 +1\right] +i \right).$$
Using this notation, we see that $\lim_{n \rightarrow \infty} \mu^{(n)} = \mu^\ast$ is a left eigenvector of the eigenvalue equal to $1$, i.e. $\mu^\ast {\bf M} = \mu^\ast$, and therefore it is invariant.  Moreover,
$$\pi_\infty (i,j;\alpha) = \mu^\ast \left(j\left[ {\rm C}^N_2 +1\right] +i \right).$$ 
\\
\\
\section*{Appendix 4: Sources of Deviation from Cellular Automaton Approximation}
This Appendix offers some detail about the deviations between the deterministic cellular automaton approximation and the fully stochastic propagation of the Markov chain paths.  On the one hand, the ergodicity of $\pi_\infty (i,j;3)$ for $N = 10$ which we just showed implies that for any function $F: {\mathcal S} \times {\mathcal A} \rightarrow {\Bbb R}$, 
$$\lim_{n \rightarrow \infty} \sum_{{\mathcal S} \times {\mathcal A}} F(x,y) \pi_n (x,y;3) = \lim_{n \rightarrow \infty} {\frac {1}{n}} \sum_{k = 1}^n F\left( S^+_k, A^+_k \right).$$
On the other hand, if we propagate (\ref{eq:xdeterm}) and (\ref{eq:ydeterm}) to create the cellular automaton $(x_n,y_n)$ we described in section 4, we evidently obtain
$$\lim_{n \rightarrow \infty} {\frac {1}{n}} \sum_{k = 1}^n F(x_n,y_n) \ne \lim_{n \rightarrow \infty} \sum_{{\mathcal S} \times {\mathcal A}} F(x,y) \pi_n (x,y;3).$$
There are two reasons for this:
\begin{enumerate}
\item \underline{Propagation effects:} 
\begin{eqnarray*}
{\bf E} \left[ \left. S^+_{k+2} \right| {\bf X}_k = (x,y) \right] & = & {\bf E} \left[ \left. {\bf E} \left[ \left. S^+_{k+2} \right| {\bf X}_{k+1} = (u,v) \right] \right| {\bf X}_k = (x,y) \right] \\
& = & \sum_{{\mathcal S} \times {\mathcal A}} {\bf E} \left[ \left. S^+_{k+2} \right| {\bf X}_{k+1} = (u,v) \right] {\rm Pr} \left( \left. {\bf X}_{k+1} = (u,v) \right| {\bf X}_k = (x,y) \right] \\
& = & {\frac {2(N-x)}{N(N+1)}} \left[ 1- P_{--}(x,y) \right] \left\{{\frac {2f(x+1,y)}{N+1}}+x+1 \right\} + \\
& & + {\frac {2x}{N(N+1)}} \left[ 1- P_{++}(x,y) \right] \left\{{\frac {2f(x-1,y)}{N+1}}+x-1 \right\} + \\
& & + {\frac {N^2-N-2y}{N(N+1)}} \left[ 1- Q_{--}(x,y) \right] \left\{{\frac {2f(x,y+1)}{N+1}}+x \right\} + \\
& & + {\frac {2y}{N(N+1)}} \left[ 1- Q_{++}(x,y) \right] \left\{{\frac {2f(x,y-1)}{N+1}}+x \right\} + \\
& & + {\frac {2L(x,y)}{N(N+1)}} \left\{{\frac {2f(x,y)}{N+1}}+x \right\}
\end{eqnarray*}
$$
\ne {\frac {2}{N+1}} f \left(\left[{\frac {2f(x,y)}{N+1}}+x \right],\left[{\frac {(N-1)g(x,y)}{N+1}}+y \right] \right) + {\frac {2f(x,y)}{N+1}}+x.$$
\item \underline{Discretization effects:} 
$$f \left(\left[(x,y)+x \right],\left[g(x,y)+y \right] \right) \ne f \left(x+{\rm sgn} \left\{f(x,y) \right\},y+{\rm sgn} \left\{g(x,y) \right\} \right),$$
\end{enumerate}
where $[x]$ denotes the nearest integer to $x$.

\end{document}